\documentclass{amsart}
\usepackage[square,sort,comma,numbers]{natbib}
\usepackage{multirow}
\usepackage{amsfonts}
\usepackage{graphicx}
\usepackage{amsmath}
\usepackage{amssymb}
\usepackage[matrix,arrow,curve]{xy}
\usepackage[colorlinks=true, allcolors=blue]{hyperref}
\usepackage[usenames]{color}
\usepackage{xcolor}
\usepackage[mathmode, boxsize=3.7em]{ytableau}
\usepackage{ulem}
\usepackage{xcolor}

\newtheorem{thm}{theorem}[section]
\newtheorem{theorem}[thm]{Theorem}
\newtheorem{proposition}[thm]{Proposition}
\newtheorem{lemma}[thm]{Lemma}
\newtheorem{corollary}[thm]{Corollary}
\newtheorem{remark}[thm]{Remark}

\newtheorem{definition}[thm]{Definition}

\newcommand{\FF}{K\{ X\}}
\newcommand{\fx}{f(x_1,\ldots,x_n)}
\newcommand{\Spec}{\mathrm{Spec}}
\newcommand{\SSpec}{\mathbb{S}\mathrm{pec}}
\newcommand{\Aschem}{\mathbb{A}\mathfrak{S}chem}
\newcommand{\Hom}{\mathrm{hom}}

\begin{document}
	
	\title[Embedding theorems: supertraces and supergeometry]{Embedding theorems as a bridge between supertraces and supergeometry}
	
	\author[Almeida]{Charles Almeida}
	\thanks{{\bf Corresponding Author: C. Fideles.}}
	\thanks{C. Almeida was supported by FAPEMIG-MG Grant No APQ-01655-22 }
	\address{Department of Mathematics, ICEx - UFMG,  Av. Ant\^onio Carlos, 6627, 30123-970 Belo Horizonte, MG, Brazil} 
	\email{charlesalmeida@mat.ufmg.br}
	
	\author[Centrone]{Lucio Centrone}
	\thanks{L. Centrone was supported by PNRR MUR project PE0000023-NQSTI}
	\address{Dipartimento di Matematica, Università degli Studi di Bari Aldo Moro, via Orabona 4, 70125, Bari, Italy}
	\email{lucio.centrone@uniba.it}
	
	\author[Fidelis]{Claudemir Fideles}
	\thanks{C. Fideles was partially supported by CNPq grant No.~305651/2021-8 and  No.~404851/2021-5 and by FAPESP grant No. 2023/04011-9}
	\address{Department of Mathematics, Universidade Estadual de Campinas (UNICAMP), 13083-859 Campinas, SP,  Brazil}
	\email{fideles@unicamp.br}

	\maketitle
	
	\begin{abstract} Any algebra herein is intended over a field of characteristic 0. Let $E$ denote the infinite dimensional Grassman algebra. Given a power associative finite dimensional {$\mathbb{Z}_2$-graded-central-simple} $A$ and a supertrace algebra $B$, so that $B$ belongs to the same variety of $A\otimes E$, we study conditions on $B$ so that it can be embedded into $A\otimes\Xi$, where $\Xi$ is a supercommutative algebra, called $A$-universal supermap of $B$, provided $B$ satisfies all the supertrace identities of $A\otimes E$. We use this result in order to relate the formal smoothness of $B$ with that of its $A$-universal supermap.

		\medskip
		
		\noindent
		\textbf{Keywords} Embedding of algebras, (super)trace identities, formal smoothness
		
		\medskip
		
		\noindent
		\textbf{Mathematics Subject Classification 2020: 16R30, 16W55, 17A05, 17A50, 14L15}
	\end{abstract}

	\section{Introduction}
	
	The development of the theory of algebras satisfying polynomial identities (PI-algebras for short) was headed by the  following  general problem: what can one say about the structure of an algebra knowing it satisfies a polynomial identity? On this purpose, many important results have been obtained throughout the last fifty years and some of them, jointly with a general overlook concerning PI-algebras, can be found in the books \cite{AGPR, KBKR, Drensky, Bahturin, varios}.
	
	A fruitful way to answer to the previous question is dealing with the following more specific and classical problem (\text{\it embedding problem}):
	\[\textit{Find necessary and sufficient conditions for a ring $R$ in order to be embedded in}\]\[\textit{a matrix algebra $M_n(C)$, where $C$ is a commutative ring.}\]
	
	A first look inside the problem shows an immediate necessary condition for embedding a ring $R$ into the algebra of $n\times n$ matrices is that it must satisfy all polynomial identities of the algebra of $n\times n$ matrices, but it is not sufficient.  For instance, we can find examples of rings that could not be embedded in a matrix algebra (and related problems) in \cite{Amitsur, Small}. Thus, a satisfactory answer to the embedding problem in general is not known yet. On the other hand, if an algebra $R$ satifies the Cayley-Hamilton identity of degree $n$, over a field of characteristic zero, we have that a ring can be embedded in a matrix algebra $M_n(C)$, provided that $C$ is a commutative ring (see \cite{procesi}). Here we recall that the Cayley-Hamilton polynomial of a matrix $a$ can be written as a polynomial whose coefficients are polynomials in the traces of $a^k$, $k\ge 1$.

	It is natural to study the embedding problem in other varieties of algebras,
	having a “richer” signature, as it was done by Berele in \cite{berele}, by Fideles, Diniz and Koshlukov in \cite{CDP} and, recently, by Almeida, Fideles and Galdino in \cite{CCGIsrael}. Surprisingly, in all of them, it was shown that if $A$ and $B$ belong to the same variety of algebras with trace, under certain additional conditions, we can embed $B$ into $A$ seen as an algebra over $C$, where $C$ is a commutative algebra so that $B$ coincides with the algebra of invariants of $A\otimes C$ by the automorphism group of $A$.

	Chapter 3.4 of the book \cite{AGPR} deals with the existence of general maps from rings to matrix algebras showing how deep the relation between noncommutative algebraic
	geometry and the theory of algebras with polynomial identities (PI-algebras) is;
	such a relation was strengthened in [4, 37]. We want to mention that the methods used to develop the solution to the embedding problem given so far relies on general maps. Such maps were introduced and studied on matrices by Amitsur, see \cite{Amitsur1}.

	{The main goal of this paper is studying embedding theorems for algebras (not necessarily associative) over an \textit{algebra of superscalars}. Here an algebra of superscalars is an unital associative which is supercommutative. The quintessential example of an algebra of superscalars is the infinite dimensional Grassmann algebra $E$. 
		{It is worth mentioning the Grassmann algebra and its canonical $\mathbb{Z}_2$-grading $E=E^{(0)}\oplus E^{(1)}$ played a key role in Kemer's theory. In \cite{basekemer2,kemer} Kemer developed a deep and far-reaching theory of varieties of associative algebras with polynomial identities where $E$ is the most powerful tool (and, at the moment, the only one) to construct a general $T$-ideal from the one of a finite dimensional algebra over a field of characteristic 0 via the so-called \textit{Grassmann envelope}. Consequently, gradings of the Grassmann algebra may shed light on what kind of mathematical construction should be investigated in order to find out how to cover identities from weaker ones. When the field is infinite and of positive characteristic, the Grassmann algebra turned out to be crucial too (see, for example, \cite{basekemer4,basekemer5,basekemer6}).
		}
		

		
	}
	
	\
	
	In this paper we prove the following (see Thorem \ref{main.theor} in the text) in the context of supercommutative algebras.

	\begin{theorem}[Main result]
		Let $A$ be a power associative algebra. Assume that $A$ is finite dimensional and central simple satisfying some additional conditions. If $B$ is an algebra with supertrace belonging to the same variety of $A\otimes E$ and satisfying all supertrace identities of $A\otimes E$, then there exists a supercommutative algebra $F_S$ such that {$B$} can be embedded into $A\otimes F_S$ as a {$\mathbb{Z}_2$-graded }algebra.
	\end{theorem}
	
	As a consequence of the previous result, we shed light on the bridge between superalgebraic geometry and PI-theory through the idea of a \textit{comm-formally smooth superscheme}. We recall a morphism of superschemes $f\colon X\to Y$ is comm-formally smooth if for
	every affine $Y$-superscheme $\SSpec(R)$ (here $R$ denotes a supercommutative algebra) and every nilpotent ideal $N \subseteq R$, any $Y$-morphism $\SSpec (R/N) \to X$ extends to an $Y$-morphism 	$\SSpec (R) \to X$. 
	
	Now, in the hypotheses of our main result, we have (see Corollary \ref{explicitcriteria}):
	\[\text{\it $B$ is formally smooth if and only if $\SSpec(F_S)$ is smooth as a superscheme.}\]
	
	Moreover, we can also relate the growth of $B$ to that of $\SSpec(F_S)$ that is the content of Corollary \ref{lastcor}: \[\text{\rm GK}(B)\leq\dim\SSpec(F_S),\]where $\text{\rm GK(B)}$ denotes the Gelfand-Kirillov dimension of $B$ whereas $\dim\SSpec(F_S)$ stands for the classical Krull dimension of $\SSpec(F_S)$.

	{Supergeometry relies on the grading by the cyclic group of order 2, denoted by $\mathbb{Z}_2$, which is ubiquitous in mathematics. It incorporates the Grassmann algebra, defined using signs, thereby excluding characteristic 2 due to its complications in forming a comprehensive supertheory. This theory has been extensively studied in recent years. The interested reader can consult the papers \cite{BHP,DorBelovVishne} and their references for a detailed description of the subject.}
	
	Associative superalgebras and their representations, known as supermodules, offer an algebraic framework for expressing supersymmetry in theoretical physics. Although they initially appeared as graded algebras in the context of algebraic topology and homology, they have gained significance in supergeometry. In this field, they are used in defining graded manifolds, supermanifolds, and superschemes.
	
	{When considering supergeometry on non-associative algebras, no specific details have been provided yet. However, the previous remarks suggest a close connection between supergeometry and the construction of generalized notions of superalgebras and superidentities. This relationship could potentially address open questions in the theory through the lens of supergeometry. 
	}
	
	{  We hope that our results concerning embedding theorems in non-associative superalgebras, referred to as $\mathbb{Z}_2$-graded algebras, might provide further insights into this area and contribute to the understanding of this emerging field in mathematics.}

	{To conclude this introductory section,} {we want to highlight the fact we used general machinery from the theory of reductive groups and Invariant Theory. Hence extending such results to fields of positive characteristic would require developing the theory of linearly reductive groups over such fields, a task that is far out of reach at the moment. Thus, we will focus on algebras over a field of characteristic zero.}
	
	\section{Preliminaries}\label{preliminary}
	
	Every algebra herein, unless otherwise stated, is supposed to be unitary and {over a ground field }$K$; a \textit{non-associative} algebra simply means a {not necessarily} associative algebra. 
	
	{Let $V$ be a vector space; a grading by the group $G$ on $V$ (or a $G$-grading on $V$) is a vector space decomposition $V=\oplus_{g\in G}V^{(g)}$, and $V$ is said to be $G$-graded (or simply graded when the group $G$ is clear from the context). An element $v\in V$ is called \textit{homogeneous} (in the $G$-grading) if $v\in V^{(g)}$ for some $g\in G$, if $v\neq 0$ we say that $v$ is homogeneous of degree $g$ and denote $\deg(a)=g$.  In particular, if $A$ is a non-associative algebra then a $G$-grading on $A$ is a $G$-grading  $A=\oplus_{g\in G} A^{(g)}$ on the underlying vector space such that $A^{(g)}A^{(h)}\subseteq A^{(gh)}$ for every $g,h\in G$; a subalgebra (an ideal, a subspace) $W$ of $A$ is \textit{graded homogeneous} if $W=\oplus_{g\in G} (W\cap A^{(g)})$. {We say that \( A \) is {\it graded-simple} if it has no proper graded (two-sided) ideals. Note that a graded-simple algebra is not necessarily a simple algebra in the classical sense. A classical example is the algebra \( A = M_n(K \oplus cK) \), where \( c^2 = 1 \), with a \( \mathbb{Z}_2 \)-grading given by \( A^{(0)} = M_n(K) \) and \( A^{(1)} = cM_n(K) \). It is worth remarking while this algebra is semisimple in the ungraded sense and isomorphic to \( M_n(K) \otimes_K K\mathbb{Z}_2 \), it does not have any proper graded ideal.} If $A$, $B$ are $G$-graded algebras, an algebra homomorphism $\varphi\colon A\to B$ is a \textsl{homomorphism of graded algebras (or $G$-graded homomorphism)} if $\varphi(A^{(g)})\subseteq B^{(g)}$ for all $g\in G$. In particular, this gives the notion of endomorphism of graded algebras, or $G$-endomorphism for short. It is important to mention that an \textit{embedding} is an one-to-one homomorphism. Let $A=\oplus_{g\in G}A^{(g)}$ and $B=\oplus_{h\in H}B^{(h)}$ be algebras graded by the groups $G$ and $H$, respectively. The tensor product grading on $A\otimes B$ has decomposition given by subspaces $(A\otimes B)^{(g,h)} = A^{(g)}\otimes B^{(h)}$, $g\in G$, $h\in H$. {Finally, a graded algebra $A=\oplus_{g\in G} A_g$ is a \textit{graded division algebra} if every non-zero homogeneous element is invertible.} 
		
		Now, consider a homogeneous element \( r \) in a \( G \)-graded algebra \( A \). We define the left and right translations, \( r_L \) and \( r_R \), respectively, as follows:  
		\begin{equation}\label{multiplication}
		r_R(a) = ar, \quad r_L(a) = ra \quad (\text{for every } a \in A).
		\end{equation}
		Right and left translations generate a subalgebra \( \mathcal{M}(A) \) of the associative algebra \(\text{End}(A) \). From classical theory, this subalgebra is referred to as the {\it multiplication algebra of \( A \)}. 
		
		The centroid of an algebra \( A \) is defined as \( C(A) = \text{End}_{\mathcal{M}(A)}(A) \), that is, the centralizer of the multiplication algebra \( M(A) \) in the endomorphism algebra \( \text{End}(A) \). In other words, \( C(A) \) consists of the linear transformations on \( A \) that commute with the multiplication operation of \( A \). If \( A \) is a simple algebra, then \( C(A) \) will be a field extension of the base field \( F \). We say that \( A \) is central if \( C(A) = F \). In the case \( A \) is $G$-graded by a group \( G \), the homogeneous operators on \( A \) form a subalgebra of \( \text{End}(A) \), denoted by $\text{End}^G(A)$, and this subalgebra is \( G \)-graded. The graded centroid \( C_{\text{gr}}(A) \) is the centralizer of \( M(A) \) in \( \text{End}^G(A) \). If \( A \) is {\it graded-simple}, then \( C_{\text{gr}}(A) = C(A) \), and it is a graded division algebra.

		{Let \( K \) be a field of characteristic different from 2 and 3, and let \( H \) be an abelian group. A map \( \beta \colon H \times H \to K^\times \) is called a {\it bicharacter} if the following conditions hold for all \( g_1, g_2, h_1, h_2 \in H \):  
			\[
			\beta(g_1g_2, h) = \beta(g_1, h)\beta(g_2, h),
			\]
			\[
			\beta(g, h_1h_2) = \beta(g, h_1)\beta(g, h_2).
			\]
			In particular, \( \beta(\epsilon, h) = \beta(g, \epsilon) = 1 \), where \( \epsilon \) denotes the identity element of \( H \).}
		
		{A bicharacter \( \beta \colon H \times H \to K^\times \) is said to be a {\it skew-symmetric bicharacter} if, for every \( h_0 \in H \), the maps \( h \mapsto \beta(h_0, h) \) and \( h \mapsto \beta(h, h_0) \) are characters of \( H \), and \( \beta(h_2, h_1) = \beta(h_1, h_2)^{-1} \) for all \( h_1, h_2 \in H \). Finally, the bicharacter \( \beta \) is called {\it non-degenerate} if \( \beta(h_0, h) = 1 \) for all \( h \in H \) implies \( h_0 = \epsilon \).}
		
		{Let \( A = \bigoplus_{g \in G} A^{(g)} \) be a \( G \)-graded algebra. To simplify the notation, we write \( \beta(a, b) \) to denote \( \beta(\deg(a), \deg(b)) \). The $G$-graded algebra $A$ is called {\it graded-central} if $C(A)_\epsilon = K\cdot 1$, , where $\epsilon$ is the identity of $G$. Moreover, \( A \) is called {\it graded-central-simple} if it is both graded-central and graded-simple, and it is {\it graded-simple and central} if is graded-simple and central (in the ungraded sense).}
		
		We have the following definition of supercommutative algebra.
		
		\begin{definition}
			An algebra $A$ is said to be supercommutative if it is a  {$\mathbb{Z}_2$-graded algebra such that for any} homogeneous elements $a,b\in A$ we have they supercommute, i.e., $ab=(-1)^{\deg(a)\cdot\deg(b)}ba$.
		\end{definition}
		
		{For instance, an important example of a graded algebra is the Grassmann algebra. The $\mathbb{Z}_2$-grading on this algebra is widely used in various parts of Mathematics and Theoretical Physics (the so-called ``super''-structures); it played a major role in Kemer’s proof (see \cite{kemer}) of the Specht problem. The Grassmann algebra is defined as follows: let $L$ be an infinite dimensional vector space over an any field $K$ of characteristic different from two with basis $\{e_{1}, e_{2}, \dots\}$. The infinite dimensional \textit{Grassmann algebra} $E$ of $L$ over $K$ is the vector space with a basis consisting of 1 and all products $e_{i_1}e_{i_2}\cdots e_{i_k}$, where $i_{1}<i_{2}<\cdots <i_{k}$, $k\geq 1$. The \textit{length} of $e_{i_1}e_{i_2}\cdots e_{i_k}$ is the number $k$ that is denoted by $|e_{i_1}e_{i_2}\cdots e_{i_k}|$. The multiplication in $E$ is induced by $e_{i}e_{j}=-e_{j}e_{i}$ for all $i$ and $j$. We shall denote the above canonical basis of $E$ by $B_{E}$. The Grassmann algebra has a canonical $\mathbb{Z}_2$-grading $E_{can}=E^{(0)}\oplus E^{(1)}$, where $E^{(0)}$ is the vector space spanned by 1 and all products $e_{i_1}\cdots e_{i_k}$ with even $k$, while $E^{(1)}$ is the vector space spanned by the products of odd length. It is well known that $E^{(0)}=Z(E)$ (here $Z(E)$ denotes the centre of $E$), and $E^{(1)}$ is the ``anticommuting'' part of $E$. Thus, $E$ endowed with its canonical $\mathbb{Z}_2$-grading is a classic example and the best known of supercommutative algebra and, on this purpose, we remand the reader to the papers \cite{krr1,giamkos,did1,cen1} showing how supercommutativity influences the identities and the graded identities of the Grassmann algebra either in the case the field has characteristic 0 or has positive characteristic. Recall that the Grassmann algebra has different kinds of gradings by the group $\mathbb{Z}_2$ and we refer the reader to \cite{did1,agpkauto}.}

		In what follows we are going to deal with the so-called \textit{embedding problem} for power associative algebras over supercommutative algebras. {A power associative algebra $A$ over a supercomutative algebra is an algebra wich entries are superscalars such that each of its subalgebras generated by a single element is associative while $A$ is (possibly) non-associative} . We recall the embedding problem studies sufficient conditions on a ring {(or algebra)} $S$ in order to be embeddable in a fixed ring {(or algebra)} $R$.

		{We will restrict our attention to finite-dimensional \textit{central simple} algebras over a field of characteristic zero, and later we will extend the ground field to a ring of superscalar. It is important to note that the full matrix algebra \( M_n(K) \), consisting of matrices with entries from a field \( K \), is a central simple algebra for any \( n \).}

		Last but not least, we recall the main ideas concerning PI-algebras.
		
		Let $X=\{x_1,x_2,\ldots,x_n,\ldots\}$ be a (possibly countable) set of variables and denote by $\FF$ the free algebra (associative, Jordan, Lie, etc.) freely generated by the set $X$ over $K$. We refer to the elements of $\FF$ as {\it polynomials}. An ideal $I$ of $\FF$ is said to be a {\it $T$-ideal} if it is invariant under all endomorphisms of $\FF$. If $A$ is an algebra over $K$, a polynomial $\fx$ is said to be a {\it polynomial identity} of $A$ if $f(a_1,a_2,\ldots,a_n)=0$ for all $a_1,a_2,\ldots,a_n\in A$. We denote by $Id(A)$ the ideal of all polynomial identities of $A$. It is a $T$-ideal of $\FF$ in the sense that it is invariant under all homomorphism of $\FF$; we will call \textit{relatively free algebra} of $A$ the quotient algebra $\FF/Id(A)$.

		We shall call \textit{substitution} with elements of $A$ any homomorphism $\FF\rightarrow A$ and we sometimes use the notation $\overline{x}=a\in A$ in order to denote explicitly such evaluation of the variable $x$.
		
		Given a subset $S\subseteq \FF$ one can think about the least $T$-ideal of $F\langle X\rangle$ containing the set $S$. Such
		$T$-ideal will be denoted by $\langle S\rangle^{T}$ and will be called the \textit{$T$-ideal generated by $S$}. We say that elements of $\langle S\rangle^{T}$ are \textit{consequences} of elements of $S$, or simply that they follow from $S$. If $Id(A)=\langle S\rangle^{T}$ for some algebra $A$, we say that $S$ is a \textit{basis} for the polynomial identities of $A$. If the ground field $F$ is infinite, then we are allowed to consider only the multihomogeneous identities of a given algebra. If the characteristic of the ground field is 0, then we are allowed to consider only the multilinear identities of a given algebra. 
		
		{Given a $T$-ideal $I$ of $\FF$, the \textsl{variety of algebras} $\mathcal{V}$ associated to $I$ is the class
			\[\mathcal{V}=\{\text{$B$ is an algebra}\mid I\subseteq Id(B)\}.\]
			In this case, the $T$-ideal $I$ of $\mathcal{V}$ is denoted by $Id(\mathcal{V})$ and is equal to $\cap_{B\in\mathcal{V}}Id(B)$. Finally, we say the variety $\mathcal{V}$ is generated by the algebra $A$ if $Id(\mathcal{V})=Id(A)$. In the sequel we shall denote by $K^{\mathcal{V}}\{X\}$ the relatively free algebra of $\mathcal{V}$ (or free algebra in $\mathcal{V}$) in the generators from $X$ (see \cite{vas} for more details as well as \cite[Theorem 2]{varios}). Indeed, two free algebras in $\mathcal{V}$ are isomorphic if and only if their free generating sets have the same cardinality. }

		Of course, the definitions above can be repeated almost verbatim in order to define weaker kinds of identities such as graded identities, identities with traces and so on. In particular, we shall denote the free $\mathbb{Z}_2$-graded algebra freely generated over $K$ by $X\cup Y$ by $K\{X;Y\}$. We agree the variables of degree $0$ are from the set $X$ whereas those of degree $1$ are from the set $Y$. Moreover, if $A$ is any $K$-algebra, we set $A\{X\}:=A\otimes_KK\{X\}$ that is the free algebra freely generated by the set $X$ over $A$.}
	
	\section{Universal supercommutative maps}\label{universal}
	
	In this section we will handle the embedding problems for finite dimensional central simple algebras over a supercommutative algebra.

	Before going on with the details of the proof, we will outline the technique used in \cite[Section 3]{CCGIsrael} in order to study the embedding problems for finite dimensional central simple algebras graded by arbitrary groups. We highlight the previous result is characteristic-free. {In this section, the symbol $K$ denotes any field of characteristic different from two}.
	
	We recall an algebra is called \textsl{commutative-associative} if it is commutative in the variety of associative algebras. Similarly, we say an algebra is  \textsl{supercommutative-associative} if it is supercommutative in the variety of associative algebras. {in other words, a commutative-associative algebra is a comutative ring, and a superring is a supercommutative-associative algebra.}
	
	Let $\mathcal{V}$ be a variety of $G$-graded algebras.Assume that $A$ is a finite dimensional $G$-graded-simple and central algebra in $\mathcal{V}$. Following word for word what was done in \cite[Section 3]{CCGIsrael}, it is possible to show that given a\textcolor{red}{n $\mathfrak{R}$-}\textcolor{blue}{$G$-graded} algebra $R$ in $\mathcal{V}$ there is a $G$-homomorphism 
		$$j_R\colon R\to A\otimes_K S_R,$$
		where $S_R$ is a commutative-associative \textcolor{red}{$\mathfrak{R}$-}algebra satisfying the following universal property: for any $G$-homomorphism $\varsigma\colon R\to A\otimes_K F$, where $F$ is a commutative-associative \textcolor{red}{$\mathfrak{R}$-}algebra, there exists an unique homomorphism $\eta_A\colon S_R\to F$ so that the diagram 
		\[
		\xymatrix{ R \ar[d]_{j_R} \ar[r]^{\varsigma} & A\otimes_KF \\ A\otimes_KS_R \ar[ru]_{\eta_A=id_A\otimes \eta}&  },
		\]
		is commutative, where $id_A$ denotes the 
		identity map in $A$. The pair $(S_R,j_R)$ is called \textsl{$A$-universal $G$-map for $R$}. This pair satisfies an important condition: the \textit{entries} of the elements $j_R(r)$, for all $r\in R$, together with the element 1, generate the \textcolor{red}{$\mathfrak{R}$-}algebra $S_{{R}}$. Moreover, such a pair is unique, up to isomorphism. Herein by entries of $r\in R$ we mean the following: consider $j_R(r)=\sum_ja_j\otimes b_j$, where the $b_j$'s belong to $S_R$; then the set $\{b_j\}_j$ is the set of entries of $r$. 
		
		From now on, $\mathcal{V}$ will denote a variety of algebras {equipped with a $\mathbb{Z}_2$-grading}; we fix a finite dimensional simple algebra $A\in\mathcal{V}$ with homogeneous $K$-basis $\mathcal{B}_A=\{u_1,\ldots,u_a, v_1,\ldots, v_b\}$, where $\deg(u_i)=0$ and $\deg(v_i)=1$.
		
		Moreover, if $C$ is a supercommutative algebra over $K$, we will always consider $A\otimes C$ as a{n} algebra endowed with the $\mathbb{Z}_2$-grading given by $(A\otimes {C})^{(0)}=A^{(0)}\otimes {C}^{(0)}+A^{(1)}\otimes {C}^{(1)}$ and $(A\otimes {C})^{(1)}=A^{(0)}\otimes {C}^{(1)}+A^{(1)}\otimes {C}^{(0)}$.

		The \textit{supercenter} of a $\mathbb{Z}_2$-graded algebra $A$, denoted by $\mathfrak{Z}=\mathfrak{Z}(A)$, is the set of homogeneous elements of $A$ which supercommmute with every homogeneous element of $A$, i.e., a homogeneous $a$ will be in $\mathfrak{Z}$ if and only if $ab=(-1)^{\deg a\cdot\deg b}ba$ for all homogeneous $b\in A$. As the field $K$ lies in $\mathfrak{Z}$, we have $\mathfrak{Z}$ is non-empty. Moreover, $\mathfrak{Z}$ is a supercommutative $\mathbb{Z}_2$-graded subalgebra of $A$. 
		
		\begin{definition}
			
			Let $R$ and $S$ be superrings and let $\eta\colon R\to S$ be a $\mathbb{Z}_2$-graded homomorphism. Denote by $\eta_A$ the $\mathbb{Z}_2$-homomorphism $(Id_A\otimes \eta)\colon A\otimes R  \to A\otimes S$; we say $\eta_A$ is the map induced by $\eta$ on $A$.
		\end{definition}
		

		The next is the main result in this section. We recall herein the word \textit{graded} always means $\mathbb{Z}_2$-graded.
		
		\begin{theorem}\label{Univsup}
			Let $\mathcal{V}$ be a variety of graded algebra. Assume that $A$ is a finite-dimensional graded-simple and central algebra such that $A\otimes E$ lies in  	$\mathcal{V}$. For all $B\in \mathcal{V}$,  there exists a superring $\Xi$, and a graded homomorphism $\varrho_\Xi\colon B\rightarrow A\otimes \Xi$ such that: 	
			\begin{enumerate}
				\item [(a)] \label{rho-asup} The entries of the elements of the set $\{\varrho_\Xi(b)\mid b\in B\}$ do generate, together with the element 1, the algebra $\Xi$;
				\item [(b)] For every supercommutative-associative algebra $F_S$ over $K$ and for every graded homomorphism $\kappa\colon B\rightarrow A\otimes F_S$ there exists a graded homomorphism $\zeta\colon \Xi\to F_S$ such that, for the induced map $\zeta_A\colon A\otimes \Xi\to A\otimes F_S$, the relation $\zeta_A\varrho_\Xi=\kappa$ holds.
			\end{enumerate}
		\end{theorem}
		
		{\begin{remark}
				We draw the reader's attention to the fact we are not requiring $A$ and $B$ lie in the same variety. Of course, this is obvious if we consider associative algebras, but it is not true in a general setting. For example, if $A=B_n$ (a Jordan algebra of a non-degenerate symmetric bilinear form on an $n$-dimensional vector space), we have $A\otimes E$ is not commutative, hence it is not a Jordan algebra, then $A$ lies in the variety of Jordan algebra but $A\otimes E$. Furthermore, the graded homomorphism $\varrho_\Xi$ given in the previous theorem, may be zero. 
			\end{remark}}
			
			{\begin{remark}
					Consider a variety of algebras $\mathcal{V}$, then we say a
					graded algebra $A = A^{(0)} \oplus A^{(1)}$ is a $\mathcal{V}$-superalgebra if 
					$E(A):=A^{(0)}\otimes E^{(0)} + A^{(1)}\otimes E^{(1)}\in \mathcal{V}$. The algebra $E(A)$ is called Grassmann
					envelope of $A$. Assume now $A\otimes E$ lies in a variety $\mathcal{V}$, then it is easy to see $E(A)$ also lies in $\mathcal{V}$, since $E(A)$ is a subalgebra of $A\otimes E$. From the last statement it follows $A$ is a $\mathcal{V}$-superalgebra. We alert the reader some authors use the word superalgebra to define $\mathbb{Z}_2$-graded algebras and, of course, the two definitions are related. The interested reader can consult \cite[Chapter 1]{booknonassocitive} for further details about the last comment.
				\end{remark}
			}	
			\begin{definition}\label{defisuper}
				Any pair $(\Xi,\varrho_\Xi)$ as in the proof of Theorem \ref{Univsup} will be called an $A$-universal supermap for $B$. 
			\end{definition}
			
			{Sometimes}, the algebra $\Xi$ will be denoted by $S_A(B)$, and we shall refer to the map $\varrho_\Xi$ by $\varrho_B$. Moreover, any $\zeta$ satisfying item $\text{\rm (b)}$ of Theorem \ref{Univsup} will said to be \textit{compatible} with the $A$-universal supermap $\varrho_\Xi$.
			
			\
			
			{Before deducing the existence of the $A$-universal supermap $(S_A(B),\varrho_B)$, we shall establish the uniqueness of this pair. The next result fulfills this purpose.}	
			\begin{proposition}\label{universalunicidadesup}
				Let $(\Xi,\varrho)$ and $(\Xi^\prime,\varrho^\prime)$ be $A$-universal supermaps for $B$. Then there exists a graded isomorphism $\varsigma\colon \Xi\rightarrow \Xi^\prime$ such that $\varrho^\prime=\varsigma_A\varrho$. 
				Moreover, the map $\zeta$ from item $(b)$ of Theorem \ref{Univsup} is uniquely determined for any fixed $(\Xi,\varrho)$ and $\kappa$. In particular, the $A$-universal supermap $(\Xi, \varrho)$ is unique, up to a graded isomorphism in the sense described above. 
			\end{proposition}
			\begin{proof}
				First we prove the uniqueness of the map $\zeta$. Let $\kappa\colon B\rightarrow A\otimes F_S$ be a {graded} homomorphism. Suppose there exist graded homomorphisms $\zeta\colon \Xi\rightarrow F_S$ and $\zeta^\prime\colon \Xi\rightarrow F_S$ such that ${\zeta}_{A}\varrho_\Xi={\zeta}_{A}^{\prime}\varrho_\Xi=\kappa$. Then ${\zeta}_{A}\varrho_\Xi(b)={\zeta}_{A}^{\prime}\varrho_\Xi(b)$ for every $b\in B$, therefore $\zeta(s)=\zeta^{\prime}(s)$ whenever $s$ is an entry of $\varrho_\Xi(b)$. We conclude $\zeta=\zeta^\prime$, since the set of entries of the elements in $\{\varrho_\Xi(b)\mid b\in B\}$ generates the algebra $\Xi$.
				
				Now let $(\Xi,\varrho)$ and $(\Xi^\prime,\varrho^\prime)$ be two $A$-universal supermaps for $B$. Then there exist {graded} homomorphisms $\varsigma\colon \Xi\rightarrow \Xi^\prime$ and $\varsigma^\prime\colon  \Xi^\prime\rightarrow \Xi$ such that $\varsigma_A\varrho=\varrho^\prime$ and ${\varsigma}_{A}^{\prime}\varrho^\prime=\varrho$. Hence $$({\varsigma}^{\prime}{\varsigma})_{A}\varrho={\varsigma}_{A}^{\prime}{\varsigma}_{A}\varrho=\varrho=(Id_\Xi)_A\varrho,$$ therefore it follows from the first part of the proof that $\varsigma^\prime\varsigma=Id_\Xi$. Similarly, $\varsigma\varsigma^\prime=Id_{\Xi^\prime}$ and hence $\varsigma$ is an isomorphism (with inverse $\varsigma^{\prime}$). 
			\end{proof}
			
			{We carry on through several technical lemmas. We recall we shall denote by $[a,b]$ the Lie product product $ab-ba$, whereas $a\circ b$ will denote the Jordan product $ab+ba$.} {It is worth emphasizing that, although the proofs of the following chain of lemmas are similar to those in \cite[Lemma 3.6 and Lemma 3.8]{CCGIsrael}, we decided to perform all the computations explicitly due to key modifications in the super setting that require attention.} 
			
			{\begin{lemma}\label{universalfree}
					Let $\mathcal{V}$ be a variety of graded algebras. Assume that $A$ is a finite-dimensional {graded-simple and central algebra} such that $A\otimes E$ lies in  	$\mathcal{V}$. Then there exists the $A$-universal supermap for each free graded algebra $K\{X; Y\}$ in $\mathcal{V}$.
				\end{lemma}
				\begin{proof}
					Let $\dim A^{(0)}=a$ while $\dim A^{(1)}=b$. Consider the following sets of variables:
					$$S =\{s_{i,j}\mid i>0 \textit{ and } 1\leq j \leq a+b\}$$
					and 
					$$T =\{t_{i,j}\mid i>0 \textit{ and } 1\leq j \leq a+b\}.$$
					We denote by $K\langle S\cup T\rangle$ the free associative algebra freely generated by the set $S\cup T$ over $K$, and we define a $\mathbb{Z}_2$-grading on $K\langle S\cup T\rangle$ setting $\deg s_{i,j}=(0)$ and $\deg t_{k,l}=(1)$. Let $P$ be the ideal of $K\langle S\cup T\rangle$ {generated by}
					\[[s_{i,j}, s_{k,l}]~~,~~~~[s_{i,j}, t_{k,l}],\   t_{i,j}\circ t_{k,l},\]
					for every $i$, $j$, $k$, $l$. Define $S_A(K\{X;Y\})= K\langle S \cup T \rangle/P$. We shall use the same letters $s_{i,j}$ and $t_{k,l}$ for the images of $s_{i,j}$ and $t_{k,l}$ under the projection $K\langle S \cup T\rangle\to S_A(K\{X;Y\})$. It follows that $S_A(K\{X;Y\})$ is a graded (associative) algebra. Furthermore, it is clear that $S_A(K\{X;Y\})$ is the free supercommutative algebra; its even variables are the $s_{i,j}$'s while the $t_{k,l}$'s are odd variables. 
					Notice that the set $S$ generates a free commutative-associative algebra $K[S]$ while $T$ generates the Grassmann algebra $E$, and the elements of $S$ commute with those of $T$. In other words, the free supercommutative algebra $S_A(K\{X;Y\})$ is isomorphic to $ K[S]\otimes E$.
					\newline
					Now we take $\mathcal{B}_A=\{u_1,\ldots,u_a, v_1,\ldots, v_b\}$ a basis of the vector space $A$ over $K$ such that $\deg(u_i)=(0)$ and $\deg(v_i)=(1)$, and we define $\mathcal{E}$ as the {algebra generated} by
					all generic elements of $A$ over $S_A(K\{X;Y\})$
					\begin{equation}\label{elementgen0}
					X_j=u_1\otimes s_{1,j}+\cdots+u_{a}\otimes s_{a,j}+v_1\otimes t_{(a+1),j}+\cdots+v_{b}\otimes t_{(a+b),j}
					\end{equation}
					and
					\begin{equation}\label{elementgen1}
					Y_j=u_1\otimes t_{1,j}+\cdots+u_{a}\otimes t_{a,j}+v_1\otimes s_{(a+1),j}+\cdots+v_{b}\otimes s_{(a+b),j}
					\end{equation}
					In this case, $\mathcal{E}$ is a graded subalgebra of $A\otimes S_A(K\{X;Y\})$. We consider the algebra homomorphism  $\psi_{K\{X;Y\}}$  defined by the composition map
					\begin{equation}\label{eq.1}
					K\{X;Y\}\rightarrow \mathcal{E}\hookrightarrow A\otimes S_A(K\{X;Y\}),
					\end{equation}
					given by $\psi_{K\{X;Y\}}(x_i)=X_i$ and $\psi_{K\{X;Y\}}(y_i)=Y_i$. Notice that $\psi_{K\{X;Y\}}$ is a graded homomorphism. It is important to mention that if $X\cup Y$ is an uncountable set, we would get  $\psi_{K\{X;Y\}}(x_\alpha)\in \mathcal{E}^{(0)}$ as well as $\psi_{K\{X;Y\}}(y_\alpha)\in \mathcal{E}^{(1)}$, if $\alpha\in\Gamma_g\setminus \mathbb{N}$. It is immediate that the pair $(S_A(K\{X;Y\}),\psi_{K\{X;Y\}})$ satisfies condition (a) of Theorem \ref{Univsup}. It remains to prove that condition (b) also holds. 
					
					Since $S_A(K\{X;Y\})$ is the free supercommutative--associative algebra freely generated by the elements $S\cup T$ it follows that for every supercommutative--associative algebra $F$ over $K$ and every $\sigma\colon K\{X;Y\}\to A\otimes F$ given by $\sigma(x_j)=\sum (u_i\otimes f_{i,j}^{(0)})+\sum (v_k\otimes f_{a+k,j}^{(1)})$ and $\sigma(y_j)=\sum (u_i\otimes f_{i,j}^{(1)})+\sum (v_k\otimes f_{a+k,j}^{(0)})$, there exists a homomorphism $\eta\colon S_A(K\{X;Y\})\to F$ such that $\eta(s_{i,j}) =f_{i,j}^{(0)}$ and $\eta(t_{i,j}) =f_{i,j}^{(1)}$. The induced map $\eta_A\colon A\otimes S_A(K\{X;Y\})\to A\otimes F$ is given by 
					\begin{eqnarray*}
						\eta_A\psi_{K\{X;Y\}}(x_j)&=&\eta_A(X_j)\\	
						&=& \eta_A(u_1\otimes s_{1,j}+\cdots+u_{a}\otimes s_{a,j}+v_1\otimes t_{(a+1),j}+\cdots+v_{b}\otimes t_{(a+b),j})\\
						&=& u_1\otimes \eta(s_{1,j})+\cdots+u_{a}\otimes \eta(s_{a,j})+v_1\otimes \eta(t_{(a+1),j})+\cdots+v_{b}\otimes \eta(t_{(a+b),j})\\
						&=&\sum (u_i\otimes f_{i,j}^{(0)})+\sum (v_k\otimes f_{a+k,j}^{(1)})=\sigma(x_j).	\end{eqnarray*}
					Similarly,  $\eta_A\psi_{K\{X;Y\}}(y_j)=\sigma(y_j)$,
					and the lemma is proved.
				\end{proof}}
				
				{
					\begin{remark} \label{idealn1}
						If $B$ is another $K$-algebra, it is well known that every ideal of $A\otimes B$ is of the form
						$A\otimes I$ for some uniquely determined ideal $I$ of $B$. Moreover, the correspondence $I\to A\otimes I$ is a lattice isomorphism from the lattice of ideals of $B$ onto the lattice of ideals of $A\otimes B$. In particular, if $B$ is a simple algebra, then $A\otimes B$ is a simple $K$-algebra. These facts are standard, see for example \cite[Remark 3.3.]{CDP}.
						
						{
							It is well-known that in the classification of gradings on simple finite-dimensional algebras, graded division algebras that are simultaneously simple play a central role. According to the graded analogues of Schur’s Lemma and the Density Theorem, any such algebra is isomorphic to \( \text{End}_D(V) \), the endomorphism algebra of a finite-dimensional graded (right) vector space \( V \) over a graded division algebra \( D \). Since \( \text{End}_D(V) \) is simple, it follows that \( D \) must also be simple. Furthermore, a finite dimensional $G$-graded division algebra is isomorphic to a twisted group algebra. These results can be found, for instance, in \cite[Chapter 2]{EK}.}
						
						{
							In our specific case, \( A \) is an \( \mathcal{M}(A) \)-module, with the action defined as in \eqref{multiplication}. Furthermore, the graded submodules of \( A \) correspond precisely to its graded ideals. Assuming \( A \) is a graded-simple algebra, we deduce that \( A \) is an irreducible graded \( \mathcal{M}(A) \)-module. Consequently, \( \mathcal{M}(A) \) is isomorphic to \( \text{End}_{C(A)}(A) \). Since \( \mathcal{M}(A) \) is graded-simple (see \cite[Proposition 7]{DiogoAngel} for details), \( C(A) \) is isomorphic to a twisted group algebra \( K^\sigma H \), where \( H \) is an abelian group. The skew-symmetric bicharacter \( \beta(u, v) = \frac{\sigma(u, v)}{\sigma(v, u)} \) for all \( u, v \in H \) arises naturally in \( C(A) \). By following these steps, the previously established correspondence holds under the hypothesis of graded-simple and central algebras. The same result applies when the field \( K \) is algebraically closed, as graded-central-simple algebras coincide with graded-simple and central algebras in this setting. More generally, the conclusion holds if \( \beta \) is nondegenerate.}
					\end{remark}
				}
				{	\begin{lemma}\label{universalquociente}
						Let $A$ and $B$ be algebras satisfying the conditions of Theorem~\ref{Univsup}.
						If $(S_A(B),\varrho_B)$ is the $A$-universal supermap for $B$, 
						then for every graded ideal $I$ of $B$ there exists a graded ideal $J$ of $S_A(B)$ and a graded homomorphism  $\varrho_B^\prime\colon B/I \to A\otimes(S_A(B)/J)$ such that the pair $(S_A(B)/J,\varrho_B^\prime)$ is the $A$-universal supermap for $B/I$.
					\end{lemma}
					\begin{proof}
						Let $(\varrho_B(I))$ be the ideal of $A\otimes S_A(B)$ generated by the set $\varrho_B(I)$. Notice that $\varrho_B(I)$ is formed by homogeneous elements in $A\otimes S_A(B)$, since $\varrho_B$ is graded. The latter claim implies that $(\varrho_B(I))$ is a graded ideal. According to Remark \ref{idealn1} there exists a graded ideal $J$ of $S_A(B)$ such that $(\varrho_B(I))=A\otimes J$. Therefore the projection maps $p_I\colon B\rightarrow B/I$ and $p_J\colon A\otimes S_A(B)\rightarrow A\otimes (S_A(B)/J)$ are well defined.
						
						We denote by $\{b^{(k)}_i\mid i\in\Gamma,k\in \mathbb{Z}_2\}$ a set of generators with homogeneous elements of $B$, where $\deg b^{(k)}_i=k$. Assume that $\mathcal{B}_A=\{u_1,\ldots,u_a, v_1,\ldots, v_b\}$ is a basis of $A$ formed by homogeneous elements satisfying $\deg(u_i)=0$ and $\deg(v_i)=1$, and we consider
						$$\varrho_B(b^{(k)}_i)=u_1\otimes s^{(k)}_{1,i}+\cdots+u_{a}\otimes s^{(k)}_{a,i}+v_1\otimes s^{(k+1)}_{(a+1),i}+\cdots+v_{b}\otimes s^{(k+1)}_{(a+b),i}$$
						By condition (a) of Theorem \ref{Univsup}, the entries of the elements in $\{\varrho_\Xi(b^{(k)}_i)\mid i\in \Gamma,k\in\mathbb{Z}_2\}$, denoted by $s^{(k)}_{ij}$ and $s^{(k+1)}_{ij}$, together with 1, generate $S_A(B)$ as a supercommutative--associative algebra over $K$. Now consider the composition map 
						\[B\stackrel{\varrho_B}{\longrightarrow} (A\otimes S_A(B))\stackrel{p_J}{\longrightarrow} (A\otimes S_A(B)/J)
						\]
						given by
						$$p_J\varrho_B(b^{(k)}_i)=u_1\otimes (s^{(k)}_{1,i}+J)+\cdots+u_{a}\otimes (s^{(k)}_{a,i}+J)+v_1\otimes (s^{(k+1)}_{(a+1),i}+J)+\cdots+v_{b}\otimes(s^{(k+1)}_{(a+b),i}+J).$$
						Then $I\subseteq \ker( p_J\varrho_B)$ and, therefore,
						\begin{equation}\label{eq.2}
						\varrho_B^\prime\colon B/I\to A\otimes (S_A(B)/J)
						\end{equation}
						given by $\varrho_B^\prime(b_i^{(k)}+I)=p_J\varrho_B(b^{(k)}_i)$ is a well defined graded homomorphism. Furthermore the diagram
						\begin{equation*}
						\xymatrix{B \ar[d]_{p_I}\ar[r]^{\varrho_B} & A\otimes S_A(B)\ar[d]^{p_J}\\
							B/I\ar[r]_{\varrho_B^{\prime}} & A\otimes (S_A(B)/J)\\ }
						\end{equation*}
						is commutative. Indeed one has $\varrho_B^\prime(p_I(b_i^{(k)}))=\varrho_B^\prime(b_i^{(k)}+I)=
						p_J(\varrho_B(b_i^{(k)}))$.
						It is clear that $S_A(B)/J$ is generated, together with $1+J$, by the elements $s^{(k)}_{ij}+J$, with $i\in\Gamma$, $k\in \mathbb{Z}_2$ and $1\leq j\leq a+b$. Thus, it is sufficient to verify condition (b) of Theorem \ref{Univsup}. We fix a supercommutative--associative algebra $F$ and a graded homomorphism $\sigma\colon S_A(B)/I\to A\otimes F$. Since $(S_A(B),\varrho_B)$ is the $A$-universal supermap for $B$, we obtain that there exists a homomorphism $\eta\colon S_A(B)\to F$ such that for the induced map $\eta_A\colon A\otimes S_A(B)\to A\otimes F$, the equality $\eta_A\varrho_B=\sigma p_I$ holds. Now we shall  prove that there exists  $\bar{\eta}\colon S_A(B)/J\to F$  such that the diagram
						\begin{equation*}
						\xymatrix{B \ar[d]_{p_I}\ar[r]^{\varrho_B} & A\otimes S_A(B)\ar[d]^{p_J}\ar@/^1cm/[dd]^{\eta_A} \\
							B/I \ar[dr]_{\sigma}\ar[r]^{\varrho_B ^\prime} & A\otimes (S_A(B)/J)\ar@{-->}[d]^{\bar{\eta}_A}\\
							& A\otimes F\\ }
						\end{equation*}
						is commutative. To this end, it suffices to show that $J\subseteq \ker\eta$. We observe that $p_I(I)=0$, hence $\eta_A\varrho_B(I)=\sigma p_I(I)=0$, that is $A\otimes J\subseteq \ker\eta_A$. Then $J\subseteq \ker\eta$, and there exists a homomorphism $\bar{\eta}\colon S_A(B)/J\rightarrow F$ such that $\bar{\eta}(s+J)=\eta(s)$. Consequently the  map $\bar{\eta}_A\colon A\otimes (S_A(B)/J)\rightarrow A\otimes F$ induced by $\bar{\eta}$ satisfies $\eta_A=\bar{\eta}_Ap_J$. In order to conclude the proof of the lemma we observe that the equalities
						\[
						\sigma(b^{(k)}_i+I)=\sigma(p_I(b^{(k)}_i))=\eta_A\varrho(b^{(k)}_i)=\bar{\eta}_Ap_J\rho(b^{(k)}_i)=\bar{\eta}_A\varrho_\Xi^\prime p_I(b^{(k)}_i)=\bar{\eta}_A\varrho_\Xi^\prime(b^{(k)}_i+I)
						\]
						hold for every $b^{(k)}_i \in B$. Therefore $\sigma=\bar{\eta}_A\rho^\prime$, and we are done.
					\end{proof}}
					
					{
						Now we have all the tools needed for the proof of the main theorem.
						\begin{proof}[Proof of Theorem \ref{Univsup}.]
							By applying Lemma \ref{universalfree}, $(S_A(K\{X;Y\}),\psi_{K\{X;Y\}})$ is the $A$-universal supermap for free graded algebra $K\{X; Y\}$ in $\mathcal{V}$. Let $\{b_i^{(k)}\mid i\in\Gamma_k, k\in \mathbb{Z}_2\}$ be a generating set of homogeneous elements of $B$ and consider the surjective graded homomorphism 
							\begin{equation*}\label{eq.3}
							\varphi\colon K\{X; Y\}\to B
							\end{equation*}
							defined by $\varphi(x_i)=b^{(0)}_i$ while $\varphi(y_i)=b^{(1)}_i$. In this way, we represent $B$ as $K\{X; Y\}/I$, where $I=\ker\varphi$. Let $P$ be the image of the ideal $I$ under $\psi_{K\{X;Y\}}$. Now the proof of the theorem follows from  Lemma \ref{universalquociente}.
						\end{proof}
					}
					
					{Let  $\mathcal{S}$ denote the category of {supercommutative-associative algebra} and $\mathcal{C}$ the category of algebras in the variety of $A\otimes E$. 
					}
					
					{From what we have done so far, the construction $B \to S_A(B)$ is clearly functorial in $B$, hence $S_A: \mathcal{S} \to \mathcal{C}$ induces a functor. Indeed,  for every morphism $f \colon B \to C$, we have the following commutative diagram:
						\begin{equation*}
						\xymatrix{B \ar[d]_{f}\ar[r]^{\varrho_B} & A\otimes S_A(B)\ar[d]^{(S_A(f))_A}\\
							C\ar[r]_{\varrho_C} & A\otimes S_A(C). }
						\end{equation*}
						Here the graded homomorphism $S_A(f)\colon S_A(B)\to S_A(C)$ cames from the universal supermap. Notice that the map $S_A(f)$ is uniquely determined so that if 
						\[B\stackrel{f}{\longrightarrow} C\stackrel{g}{\longrightarrow} D,\]
						we then have
						$$S_A(g\circ f)=S_A(g)\circ S_A(f).$$
						All in all, the following result follows.
						\begin{theorem}\label{adjoint}
							The functor $S_A:\mathcal{S} \to \mathcal{C}$ has a right adjoint, more precisely we have that
							$$\hom_\mathcal{C}(B, A\otimes F) = \hom_\mathcal{S}(S_A(B), F).$$
						\end{theorem}}
						
						A nice consequence of Theorem \ref{adjoint} is that, since $S_A$ admits a right adjoint, it preserves direct limits of algebras (see \cite[Theorem 6.3.1]{TL2016})

						{We will leverage on the existence of a universal map given by Theorem \ref{Univsup} in order to relate the graded algebras considered here with superschemes in the traditional sense. Let $F$ be a superring; the \textit{superspectrum} of $F$, $\SSpec(F)$, is a pair $(X, \mathcal{O}_X)$ consisting of a topological space $X$ which is the spectrum of the even part of $F$, and $\mathcal{O}_X$, a sheaf of $\mathbb{Z}_2$-graded rings on $X$, such that for each $f \in F$, we have a basic open $D(f) \subseteq X$, and one defines $\mathcal{O}_X(D(f)) = F_f$ the localization of $F$ to the multiplicative subsystem defined by $f$.  With abuse of notation, we say that $\SSpec(F)$ is the \textit{affine superscheme associated to $F$}}. 
						{Let $\mathbb{A}\mathfrak{S}chem$ denote the category of affine superschemes, and {\bf Sets} the usual category of sets. Assume that $A$ is a finite-dimensional {graded-simple and central algebra}; we fix a graded algebra $B$ belonging to the same variety of algebras of $A\otimes E$. Consider the functor \[\mathcal{F}_B^A: \mathbb{A}\mathfrak{S}chem \to {\mathbf{\mathrm{Sets}}}\] such that for each $X = \SSpec(S) \in Obj( \mathbb{A}\mathfrak{S}chem)$ we have \[\mathcal{F}_B^A(X):= \mathrm{hom}_{\mathcal{V}}(B,A\otimes S),\] where $\mathrm{hom}_{\mathcal{V}}(B,A\otimes S)$ denotes the set of homomorphisms from $B$ to $A\otimes S$ with the obvious actions on morphisms and compositions. More precisely, if $g:\SSpec(S_1) \to \SSpec(S_2)$ is a morphism of superschemes, it induces a morphism of supercommutative algebras, $\tilde{g}: S_2 \to S_1$, therefore we can define $\mathcal{F}_B^A(g): \mathrm{hom}_{\mathcal{V}}(B,A\otimes S_2) \to \mathrm{hom}_{\mathcal{V}}(B,A\otimes S_1)$ as $\mathcal{F}_B^A(g)(\rho) = (I_A\otimes \tilde{g})\circ \rho$.
							Note that, given a morphism of supercommutative algebras $f: S_1 \to S_2$, there is a morphism of affine superschemes $f^{\#}: \SSpec(S_2) \to \SSpec(S_1)$ satisfying by definition $\widetilde{f^{\#}}= f$. If $g: S_3 \to S_2$ is another morphism of supercommutative algebras, then $\widetilde{f \circ g}= \tilde{g} \circ \tilde{f}$, and $(\tilde{g} \circ \tilde{f})^{\#} = f \circ g$.}
						
						We have the next.
						\begin{theorem}\label{suprep}
							With the notation as above, $\SSpec(S_A(B))$ represents the functor $\mathcal{F}_B^A$, that is, $\mathcal{F}_B^A \simeq\Hom_{\Aschem}(-,\SSpec(S_A(B)))$ 
						\end{theorem}
						\begin{proof}
							
							We are going to prove that there exists a natural transformation $\eta$ between $ \mathcal{F}_B^A$ and $\Hom_{\Aschem}(-,\Spec(S_A(B)))$. Recall that by Theorem 
							\ref{Univsup},  there exists an $A$-universal supermap $(S_A(B), \rho)$ for $B$.
							Now, let $\tau \in \mathcal{F}_B^A(X) =  \mathrm{hom}_{\mathcal{V}}(B,A\otimes S)$ for some $X = \SSpec(S)$, where $S$ is a supercommutative-associative algebra.  By Theorem \ref{Univsup}, there exists $\zeta_{\tau}: S_A(B) \to S$ compatible with the $A$-universal supermap $\rho$. This will induce a map $\zeta_{\tau}^{\#}: \SSpec(S) \to \SSpec(S_A(B))$. Hence we define $\eta_X(\tau) = \zeta_{\tau}^{\#}$. 
							Now we prove the next diagram is commutative.

							\[\xymatrix{ \mathcal{F}_B^A(S_2)  \ar[r]^-{\eta_{X_2}} \ar[d]_{\mathcal{F}_B^A(g)} & \Hom_{\Aschem}(\SSpec(S_2),\SSpec(S_A(B))) \ar[d]^{\circ g} \\
								\mathcal{F}_B^A(S_1) \ar[r]^-{\eta_{X_1}} & \Hom_{\Aschem}(\SSpec(S_1),\SSpec(S_A(B)))}. \]

							Given a morphism of superschemes $g: X_1:= \SSpec(S_1)\to \SSpec(S_2)=: X_2$,
							we have a morphism of supercommutative algebras $\tilde{g}: S_2 \to S_1$ and an induced map 
							$\Hom_{\Aschem}(\SSpec(S_2),\SSpec(S_A(B))) \to \Hom_{\Aschem}(\SSpec(S_1),\SSpec(S_A(B)))$ given by the composition with $g$, that we will denote by $\circ g$.

							Now, given $\tau: B \to A\otimes S_2$ and $g: X_1:= \SSpec(S_1)\to \SSpec(S_2)=: X_2$, note that we have the following commutative diagram:

							\[\xymatrix{ B \ar[rr]^{\tau} \ar[rrrd]^{\rho_B} && A \otimes S_2 \ar[rr]^{I \otimes \tilde{g}} & & A \otimes S_1\\
								& & & A \otimes S_A(B) \ar[ul]_{I\otimes \zeta_{\tau}}  \ar[ur]_{I\otimes \zeta_{(I\otimes\tilde{g})\circ \tau}}&}\]
							
							\noindent which gives us that $\zeta_{(I\otimes\tilde{g})\circ \tau} = \tilde{g} \circ \zeta_{\tau} $. This identity will imply that the following diagram commutes:
							
							\[\xymatrix{ \mathrm{hom}_{\mathcal{V}}(B,A\otimes S_2)  \ar[r]^-{\eta_{X_2}} \ar[d]_{\mathcal{F}_B^A(g)} & \Hom_{\Aschem}(\SSpec(S_2),\SSpec(S_A(B))) \ar[d]^{\circ g} \\
								\mathrm{hom}_{\mathcal{V}}(B,A\otimes S_1) \ar[r]^-{\eta_{X_1}} & \Hom_{\Aschem}(\SSpec(S_1),\SSpec(S_A(B)))} \]
							because for each $\tau: B \to A\otimes S_2$ we have
							
							$$\circ g(\eta_{X_2}(\tau)) = \circ g (\zeta_{\tau}^{\#})= \zeta_{\tau}^{\#}\circ g$$
							
							\noindent and
							
							\[\eta_{X_1}(\mathcal{F}_B^A(g)(\tau)) =  \eta_{X_1}((I_A\otimes \tilde{g})\circ \tau)= (\zeta_{(I\otimes\tilde{g})\circ \tau})^{\#} = (\tilde{g} \circ \zeta_{\tau})^{\#} = \zeta_{\tau}^{\#}\circ g.\]

							\noindent Thus $\eta_X$ provides a natural transformation from $\mathcal{F}_B^A$ to $\Hom_{\Aschem}(-,\SSpec(S_A(B)))$. Now we shall prove for each $X_S \in \Aschem$, with $X_S = \SSpec(S)$, the natural transformations $\eta_{X_S}: \mathcal{F}_B^A(X_S) \to \Hom_{\Aschem}(X_S,\SSpec(S_A(B)))$ is one-to-one and onto. 
							
							{\bf Part 1:} ($\eta_{X_S}$ is one-to-one)

							Indeed, let $\tau \in \mathcal{F}_B^A(X_S)$, and $\sigma \in \mathcal{F}_B^A(X_S)$, such that $\zeta_{\tau}^\# = \zeta_{\sigma}^\#$. But this would imply that $\widetilde{\zeta_{\tau}^\#} = \widetilde{\zeta_{\sigma}^\#}$, from which follows that $\zeta_{\tau} = \zeta_{\sigma}$. Now, we have  $$\tau = (I_A \otimes \zeta_{\tau})\circ \rho_B = (I_A \otimes  \circ \zeta_{\sigma})\rho_B = \sigma.$$
							
							{\bf Part 2:} ($\eta_{X_S}$ is onto)
							
							Given $f \in \Hom_{\Aschem}(X_S,\SSpec(S_A(B)))$, there is a $\tilde{f}: S_A(B) \to S$, and we can define $\tau_f \in \mathcal{F}_B^A(X_S)$ as $(I_A \otimes \tilde{f})\circ \rho_B$. Which concludes the proof. 
							
						\end{proof}

									Let us close the section with some consequences of Theorem \ref{Univsup}.
									
									\begin{corollary}\label{monosup}
										The algebra $B$ can be embedded into $A\otimes F_S$ over a supercommutative-associative algebra $F_S$, as a graded subalgebra, if and only if the graded homomorphism $\varrho_\Xi\colon B\to A\otimes \Xi$ of Theorem \ref{Univsup} is a monomorphism. 
									\end{corollary}
									\begin{proof}
										If $\varrho_\Xi$ is a monomorphism, then clearly $B$ can be embedded into $A\otimes \Xi$, as a graded subalgebra. Conversely, suppose there exists an embedding $\kappa\colon B\to A\otimes F_S$, where $F_S$ is a supercommutative-associative $K$-algebra, as a graded subalgebra. Then by Theorem \ref{Univsup} there exists a {graded} homomorphism $\zeta_A\colon \Xi\to F_S$ such that the diagram
										\[
										\xymatrix{ B \ar[d]_{\varrho_\Xi} \ar[r]^{\kappa} & A\otimes F_S \\ A\otimes \Xi \ar[ru]_{\zeta_A}&  }
										\]
										is commutative and $\kappa=\zeta_A\varrho_\Xi$. Hence, if $\kappa$ is a monomorphism, then  $\varrho_\Xi$ is a monomorphism as well.
									\end{proof}
									
									We shall see in Subsection \ref{sectioembeddin} a criterion realizing the hypotheses of Corollary \ref{monosup}.

									
									{For the next result, we will need to recall the notations established in the proof of Lemma \ref{universalfree}. We denote by $\mathcal{E}=\mathcal{E}_{S_A(K\{X;Y\})}$ the algebra generated by the elements $X_j$ and $Y_j$, for every $j>0$. Here $X_j$ and $Y_j$ are elements given in \eqref{elementgen0} and \eqref{elementgen1}, respectively. Such algebra is a graded subalgebra of $A\otimes {S_A(K\{X;Y\})}$, and it will be called $\mathbb{Z}_2$-graded algebra of \textit{generic superelements} on $A$.}
									
									\begin{corollary}\label{embsup}
										The algebra $B$ can be embedded, as a graded subalgebra, into $A\otimes F_S$ for some superring $F_S$ if and only if there exists a graded epimorphism $\varphi\colon {\mathcal{E}}\to B$ such that $(P)\cap {\mathcal{E}}=P$ where $P=\ker\varphi$ and $(P)$ is the graded ideal of ${\mathcal{E}}$ generated by $P$.
									\end{corollary}
									\begin{proof}
										Let us denote by {$(\Delta,\varrho_{\Delta})$} the $A$-universal supermap of ${K}\{X;Y\}$ and recall $B$ is isomorphic to the quotient $K\{X;Y\}/I$ for some graded ideal $I$ of $K\{X;Y\}$. Then there exists an ideal $J$ in {$\Delta$} such that $(P)=A\otimes J$, where $P=\varrho_{{\Delta}}(I)$. {Notice that $A\otimes J$ is graded because $P$ is graded.} By Corollary \ref{monosup}, $\varrho\colon K\{X;Y\}/I\to A\otimes({\Delta}/J)$ is a monomorphism, where $({\Delta}/J,\varrho)$ is the $A$-universal supermap for $B$. Indeed, we obtain a graded isomorphism $$\varrho_{0} \colon K\{X;Y\}/I \to {\mathcal{E}}/P.$$
										Let $\iota_P$ be the canonical homomorphism of ${\mathcal{E}}/P$ in $A\otimes ({\Delta}/J)=(A\otimes {\Delta})/(P)$, then we have $\varrho$ is given by the composition of the {graded} homomorphisms
										\[
										\xymatrix{ B\simeq K\{X;Y\}/I\ar[r]^{\ \ \ \ \varrho_{0}}& {\mathcal{E}}/P\ar[r]^{\iota_P}&  (A\otimes {\Delta}}/J) \\ .
										\]
										We consider now the graded epimorphism $\varphi\colon {\mathcal{E}}\to K\{X;Y\}/I$ given by the composition of the canonical epimorphism ${\mathcal{E}}\to {\mathcal{E}}/P$ with $\varrho_{0}^{-1}$; we have $P=\ker\varphi$. As $\varrho$ is injective, we conclude $\iota_P$ is injective too, therefore $(P)\cap {\mathcal{E}} = P$.
										
										Conversely, suppose there exists a graded epimorphism $\varphi$ from ${\mathcal{E}}$ to $K\{X;Y\}/I$ such that $P=\ker\varphi$ and satisfying $(P)\cap {\mathcal{E}}=P$. Then ${\mathcal{E}}/P\backsimeq K\{X;Y\}/I$,
										and $(P)\cap {\mathcal{E}} = P$ and now the proof follows.
									\end{proof}
									
									\begin{corollary}\label{mergsup}
										Every {graded} algebra $B$ contains a graded ideal $Q$ such that the quotient $B/Q$ can be embedded into $A\otimes {S_A(B)}$, as a graded subalgebra. {Additionally, if $Q_0$ is a graded ideal of $B$, such that}  $B/Q_0$ can be embedded into $A\otimes F_S$ for some superring $F_S$, as a graded subalgebra, then $Q\subseteq Q_0$.
									\end{corollary}
									\begin{proof}
										{Let $(S_A(B),\varrho_B)$ be the $A$-universal supermap of $B$,}
										and define $Q=\ker\varrho_B$. Therefore $\varrho_B$ induces a monomorphism of $B/Q$ into $A\otimes {S_A(B)}$. Now suppose that there exists another graded ideal of $B$, to say $Q_0$, such that $B/Q_0$ can be embedded into $A\otimes {S_A(B)}$, as a graded subalgebra. In other words $\kappa\colon B/Q_0\to A\otimes F_S$
										is a {graded} monomorphism for some supercommutative  algebra $F_S$. According to Theorem \ref{Univsup}, there exists $\zeta_A\colon A\otimes {S_A(B)}\to A\otimes  F_S$ such that $\zeta_A\varrho_B=\kappa p_0$ where $p_0\colon B\to B/Q_0$ stands for the canonical projection. Notice that for every $x\in Q$ we have $0=\zeta_A\varrho_B(x)=\kappa p_0(x)$, and this implies $Q\subseteq Q_0$ since $\kappa$ is a monomophism. The uniqueness of such an ideal is straightforward and we are done.
									\end{proof}
									
									{For the next application the definition of notherian comutative-associative algebra will be extended to supercomutative-associative algebra.}
									\begin{definition}
										A superring $\mathcal{R}$ is called Noetherian if its {graded} ideals satisfy the ascending chain condition (ACC).
									\end{definition}
									
									As it has been proven in \cite[Lemma 1.4]{MasZub}, a superring $\mathcal{R}=\mathcal{R}^{(0)}\oplus\mathcal{R}^{(1)}$ is Noetherian if and only if $\mathcal{R}^{(0)}$ is a Noetherian ring and the $\mathcal{R}^{(0)}$-module $\mathcal{R}^{(1)}$ is finitely generated.
									
									\begin{corollary}\label{usefulGK}
										In the hypotheses of Theorem \ref{Univsup}, if $B$ is a finitely generated graded algebra, then so is ${S_A(B)}$. More precisely, ${S_A(B)}$ is a Noetherian superring. 
									\end{corollary}
									\begin{proof}
										Assume $B$ is a finitely generated graded algebra, then condition (a) of Theorem \ref{Univsup} implies
										that ${S_A(B)}$ must be finitely generated as well. {Due to the last comment we consider that the free supercomutative-associative algebra has a finite amount of generators of degree 0 and 1. Thus, the set $S=\{s_1,\ldots,s_n\}$ contains commutative variables while $T=\{t_1,\ldots,t_n\}$ contains anticommutative variables. }The proof of Theorem \ref{Univsup} implies that ${S_A(B)}^{(0)}$ ($\simeq ${$R/J$ where $R=K[S]\otimes K[t_it_j\mid 1\leq i<j\leq m]$ and $J$ is an ideal of $R$)} is a Noetherian ring. Thus, it is enough to show that ${S_A(B)}^{(1)}$ is a Noetherian ${S_A(B)}^{(0)}$-module. As ${S_A(B)}$ is finitely generated, there exist elements $\mathfrak{r}_1$, \ldots, $\mathfrak{r}_n$ that generate ${S_A(B)}$ as a graded algebra. Assume, without loss of generality, that $\mathfrak{r}_1$, \ldots, $\mathfrak{r}_s$ are elements of ${S_A(B)}^{(0)}$ and $\mathfrak{r}_{s+1}$, \ldots, $\mathfrak{r}_n$ are ${S_A(B)}^{(1)}$. We will show that ${S_A(B)}^{(1)}$ is the ${S_A(B)}^{(0)}$-module generated by the set $\{\mathfrak{r}_{s+1}, \ldots, \mathfrak{r}_n\}$. By the definition of graded algebra ${S_A(B)}$, we have that $\Xi^{(1)}$ is generated as ${S_A(B)}^{(0)}$-module by all monomials of the form
										$$\mathfrak{r}_{s+1}^{\alpha_{s+1}}\cdots \mathfrak{r}_n^{\alpha_{n}}$$
										such that ${\alpha_{s+1}}$, \ldots, ${\alpha_{n}}$ are non-negative integers such that ${\alpha_{s+1}}+\cdots+{\alpha_{n}}$ is odd. Thus, there exists an index $i$ between $s+1$ and $n$ such that $\alpha_i\neq 0$, but $\alpha_{j}=0$ for all $i<j\leq s+1$. We have $\mathfrak{r}_{s+1}^{\alpha_{s+1}}\cdots \mathfrak{r}_i^{\alpha_{i}-1}$ has even length, i.e., it lies in ${S_A(B)}^{(0)}$, and we are done.
									\end{proof}
									
									The previous result shares light on the growth behaviour of $B$ with respect to the algebra ${S_A(B)}$, too. To show this, in the sequel we are going to give an upper bound of the so-called \textit{Gelfand-Kirillov (GK) dimension} of the algebra $B$. For more details about the GK dimension we remand to the classical books \cite{Krle} by Krause and Lenagan and \cite{MR} by McConnell and Robson. Here we simply point out the GK dimension of an algebra measures its rate of growth in terms of any generating set. It could be defined for modules and groups, too.
									
									The next result shows how close is the rate of growth of the algebra $B$ to that of the algebra ${S_A(B)}^{(0)}$.
									
									\begin{corollary}\label{GK2}
										In the hypotheses of Theorem \ref{Univsup}, let $B$ be finitely generated. Then \[\text{\rm GK}({\varrho_\Xi}(B))\leq \text{\rm GK}({S_A(B)}^{(0)}).\] 
										
									\end{corollary}
									\begin{proof}
										We observe 
										the GK dimension of ${S_A(B)}$ as a module over itself and as an algebra coincide, then we shall use the same symbol $\text{\rm GK}({S_A(B)})$. By Theorem \ref{Univsup} we have ${\varrho_B}(B)$ is a subalgebra of $A\otimes S_A(B)$, then by Propositions 8.2.2 and 8.2.3 of \cite{MR}, we have \begin{equation}\label{GK1}\text{\rm GK}(\varrho_B(B))\leq \text{\rm GK}(A)+\text{\rm GK}({S_A(B)}).\end{equation} 
										By Proposition 5.5 of \cite{Krle}, because ${S_A(B)}$ is a finitely generated ${S_A(B)}^{(0)}$-module and $A$ is finite dimensional, we have $\text{\rm GK}({S_A(B)})=\text{\rm GK}({S_A(B)}_{{S_A(B)}^{(0)}})$ and $\text{\rm GK}(A)=0$; moreover, by Proposition 8.3.2 of \cite{MR} and again by Proposition 5.5 of \cite{Krle} (because ${S_A(B)}^{(1)}$ is a finitely generated ${S_A(B)}^{(0)}$-module), we have $\text{\rm GK}({S_A(B)}_{{S_A(B)}^{(0)}})=\max\{\text{\rm GK}({S_A(B)}^{(0)}),\text{\rm GK}({S_A(B)}^{(1)})\}=\text{\rm GK}({S_A(B)}^{(0)})$, too. By Equation (\ref{GK1}), we get \[\text{\rm GK}(\varrho(B))\leq \text{\rm GK}({S_A(B)}^{(0)})\] and the proof follows.

									\end{proof}

									We want to give another application of Theorem \ref{Univsup} in the language of invariant theory.
									
									Let $\mathcal{G}=Aut^{\mathbb{Z}_2}(A)$ be the group of the graded automorphism of $A$, and let  {$(S_A(B),\varrho_B)$} be the $A$-universal supermap for some graded algebra $B$. Then one can consider the map $\Phi_g\colon A\otimes    {S_A(B)}\rightarrow A\otimes   {S_A(B)}$ defined by $\Phi_g(a\otimes \mathfrak{r})=(g\cdot a)\otimes \mathfrak{r}$ for every $g\in \mathcal{G}$. It follows immediately from the definition of the $A$-universal supermap that there exists a map $\zeta^{g}\colon   {S_A(B)}\rightarrow   {S_A(B)}$ such that
									\begin{equation}\label{eq.4}
									\Phi_g\circ \varrho_B=\zeta^{g}_A\circ \varrho_B.
									\end{equation}
									This implies that the diagram
									\begin{equation*}
									\xymatrix{B \ar[d]_{\varrho_B}\ar[r]^{\varrho_B} & A\otimes   {S_A(B)}\ar[r]^{\Phi_g} & A\otimes   {S_A(B)}\\
										A\otimes   {S_A(B)}\ar[urr]_{\zeta^{g}_A} & \\ }
									\end{equation*}
									is commutative. If $g$, $h\in \mathcal{G}$, then it holds $\Phi_{g}\circ\zeta^{h}_A=\zeta^{h}_A\circ\Phi_{g}$ since both $\Phi_{g}$ and $\zeta^{h}_A$ fix the elements of $S_A(B)$ and $A$, respectively. Moreover,
									\[
									\zeta^{gh}_A\circ \varrho=\Phi_{gh}\circ \varrho=\Phi_{g}\circ\Phi_{h}\circ \varrho=\Phi_{g}\circ\zeta^{h}_A\circ \varrho=\zeta^{h}_A\circ\Phi_{g}\circ \varrho=\zeta^{h}_A\circ\zeta^{g}_A\circ \varrho=(\zeta^{h}\zeta^{g})_A\circ \varrho
									\]
									implying $\zeta^{gh}=\zeta^{h}\zeta^{g}$. Now, setting $\xi_g=\Phi_g\circ(\zeta^{g}_A)^{-1}$, we obtain the following.
									
									\begin{proposition}\label{actionofGsup}
										Let {$(S_A(B),\varrho_B)$} be the $A$-universal supermap for $B$. Then $\mathcal{G}=Aut^{\mathbb{Z}_2}(A)$ acts, via the representation $\xi_g$, as a group of {graded} automorphisms for $A\otimes {S_A(B)}$. Moreover, $\varrho_B(B)$ is contained in the ring of invariants $(A\otimes {S_A(B)})^\mathcal{G}$.
									\end{proposition}
									
									We recall if $f:X\rightarrow Y$ is a map between two sets with a  $\mathcal{G}$-action, then $f$ is said to be \textit{$\mathcal{G}$-equivariant} if 
									$f(\mathfrak{g}\cdot x)=\mathfrak{g}\cdot f(x)$ for all $\mathfrak{g}\in \mathcal{G}$. Whenever $\mathcal{G}$ is understood from the context we may also say $f$ is an \textsl{equivariant map}.
									
									\begin{corollary}\label{equi}
										Let $({S_A(B)},\varrho_B)$ be the $A$-universal supermap for a graded algebra $B$. If $I$ is a graded ideal of $B$ such that $(\varrho_B(I))=A\otimes J$, then the projection map $p_J\colon A\otimes_K {S_A(B)}\rightarrow A\otimes ({S_A(B)}/J)$ is $\mathcal{G}$-equivariant.
									\end{corollary}
									\begin{proof}
										With abuse of notation we shall use the same symbol $\Phi_g$ to denote the graded homomorphism $\Phi_g\colon A \otimes ({S_A(B)}/J) \to A \otimes({S_A(B)}/J)$ given by $\Phi_g(a\otimes (\mathfrak{r}+J)) = (g \cdot a) \otimes (\mathfrak{r}+J)$ for all $ g \in \mathcal{G}$.  In order to prove the result it is enough to verify that  $p_J\circ\Phi_g\circ({\zeta}_{A}^{g})^{-1}=\Phi_g\circ(\bar{\zeta}_{A}^{g})^{-1}\circ p_J$ where $\bar{\zeta}_{A}^{g}$ is the graded homomorphism from $B/I$ into $A\otimes({S_A(B)}/J)$ given by $\bar{\zeta}_{A}^{g}=p_J\circ{\zeta}_{A}^{g}$. Observe that the diagram
										\[\xymatrix{ A\otimes {S_A(B)} \ar[r]^{\Phi_g} & A\otimes {S_A(B)} \\
											B \ar[d]_{p_I}\ar[u]^{\varrho} \ar[r]^{\varrho} & A\otimes {S_A(B)}\ar[u]_{{\zeta}_{A}^{g}}\ar[d]^{p_J} \\
											B/I  \ar[d]_{\varrho_0}\ar[r]_{\varrho_0}& A\otimes({S_A(B)}/J)\ar[d]^{\bar{\zeta}_{A}^{g}}\\
											A\otimes({S_A(B)}/J)  \ar[r]_{\Phi_g} & A\otimes({S_A(B)}/J) \\  }
										\]
										is commutative. Therefore 
										\begin{eqnarray*}
											p_J\circ\Phi_g\circ({\zeta}_{A}^{g})^{-1}\circ\varrho&=&p_J\circ({\zeta}_{A}^{g})^{-1}\circ \Phi_g\circ \varrho\\
											&=&p_J\circ\varrho =\varrho_0\circ p_I = (\bar{\zeta}_{A}^{g})^{-1}\circ\Phi_g\circ\varrho_0\circ p_I\\
											&=&(\bar{\zeta}_{A}^{g})^{-1}\circ\Phi_g\circ p_J\circ\rho =\Phi_g\circ(\bar{\zeta}_{A}^{g})^{-1}\circ p_J\circ\varrho.
										\end{eqnarray*}
										The result now follows from Proposition \ref{universalunicidadesup}.
									\end{proof}
									
									We can also think about $A\otimes E$ as the ring of polynomial maps $f\colon (A\otimes \mathcal{C})^k\rightarrow (A\otimes \mathcal{C})$ over $K$, where $\mathcal{C}$ is any supercommutative $K$-algebra. The group action is the natural one, i.e., $f^g(a_1,\dots,a_k):=g\cdot f((g^{-1}\cdot a_1,\dots, g^{-1}\cdot a_k)$ so that $f=f^g$ means $f$ is an \textit{equivariant map}.

									\section{Embedding Theorem on certain {graded} algebras {over supercommutative ring}}
									
									To obtain the embedding theorem, we will need some general machinery from trace algebras and invariant theory that we will present at the beginning of this section. 
									
									\subsection{Brief review on supertrace identities}\label{traceidentities} 
									
									Let $A$ be a non-associative algebra over ${K}$, and let $\mathcal{Z}=\mathcal{Z}(A)$ be its \textit{associative and commutative center}. A $K$-linear map $\tau\colon A\to \mathcal{Z}$ such that $\tau(\tau(a)b) = \tau(a)\tau(b)$ and $\tau([a, b]) = \tau((a, b, c)) = 0$, for all 
									$a$, $b$ and  $c$ in $A$, is called a \textsl{trace} on $A$. Here we recall $(a, b, c) = (ab)c-a(bc)$ stands for the \textit{associator} of $a$, $b$, $c\in A$. The pair $(A,\tau)$ will be called a \textsl{trace algebra} (or algebra with trace). The \textsl{pure trace algebra} for $(A,\tau)$ is defined as the image of $A$ under $\tau$, denoted by $\tau(A)$. It follows immediately that $\tau(A)$ is a subalgebra of $\mathcal{Z}$. An ideal (resp. subalgebra) in a trace algebra $(A,\tau)$ is called a \textsl{trace ideal} (resp. \textsl{trace subalgebra}) of $(A,\tau)$ if it is an ideal (resp. subalgebra) that is closed under the trace map.  Let $(A, \tau)$ and $(B,\theta)$ be trace algebras, a map $\varphi\colon(A, \tau)\to (B,\theta)$ is a \textsl{homomorphism of trace algebras} if it is a  homomorphism of algebras and $\varphi(\tau(a)) = \theta(\varphi(a))$ holds for every $a$ in $A$. It is
									worth noting that we do not require the trace map to be non-degenerate
									and the algebra, along with its subalgebras, may have different units.
									
									Denote by $[x, y]_s = xy-(-1)^{\deg x\cdot\deg y}yx$ the supercommutator of the homogeneous elements $x$, $y$. When at least one of the elements $x$, $y$ is \textit{even} (has degree $(0)$), we will often omit the index $s$ and we simply write $[x, y]$. A linear map $\tau\colon A \to \mathfrak{Z}(A)$, where $A$ is a {graded} algebra, is called a \textit{supertrace} (Str) if it satisfies $\tau(\tau(a)b) = \tau(a)\tau(b)$ and $\tau([a, b]_s) = \tau((a, b, c)) = 0$, for all homogeneous $a$, $b$ and  $c$ in $A$. Thus, a {graded} algebra with a supertrace will be called a \textit{supertrace algebra} (or \textit{{graded} algebra with supertrace}). {It is important to mention that a \textsl{supertrace ideal} (resp. \textsl{supertrace subalgebra}) of a supertrace algebra $(A,\tau)$ is a graded ideal (resp. graded subalgebra) of $A$ that is closed under the supertrace map.} A classical example of supertrace algebra is the algebra of $n\times n$ matrices over the Grassmann algebra $M_n(E)$, with $\mathbb{Z}_2$-grading inherited from the canonical $\mathbb{Z}_2$-grading of $E$. In this case, the supercenter is equal to $E\cdot I_n$ and the supertrace will be given by the usual way, i.e., for all $(a_{ij})^{n}_{i,j=1}\in M_n(E)$, the map \[\mathrm{Str}((a_{ij})^{n}_{i,j=1})=\sum_{i=1}^{n}a_{ii}\] is a supertrace of $M_n(E)$. {In general, } let $(A,\mathrm{tr})$ be any trace algebra, then $A\otimes E$ will be a {graded} algebra with supercenter ${\mathcal{Z}}\otimes E$ and with supertrace $\mathrm{tr}\otimes id$, i.e., $\mathrm{Str}(a\otimes e)=\mathrm{tr}(a)\otimes e${, where $a\in A$ and $e$ is a monomial in $E$. This definition extends linearly to the whole algebra.}

									Given two countable sets of variables $X$ and $Y$, where the variables in $X=\{x_{1}, x_{2},\ldots\}$ are of degree $0$, and those in $Y=\{y_{1}, y_{2},\ldots\}$ are of degree $1$, we define the {\it free supertrace algebra} in $\mathcal{V}$, denoted by ${K}^{\mathcal{V}}_{STR}\{X;Y\}$, in the following way. 
									
									\begin{definition}Let $K^{\mathcal{V}}\{X;Y\}$ be the free {graded} algebra in $\mathcal{V}$, then the free supercommutative supertrace associative algebra $STR^{\mathcal{V}}\langle X;Y\rangle$ is the algebra generated by the symbols $\mathrm{Str\,}(g)$, $g\in K^{\mathcal{V}}\{X;Y\}$ satisfying the next conditions:
										\begin{enumerate}
											\item[(i)] if $g$ is homogeneous, then  $\mathrm{Str\,}(g)$ is homogeneous of the same degree;
											\item[(ii)] $\mathrm{Str\,}(f+g)=\mathrm{Str\,}(f)+\mathrm{Str\,}(g)$ {and $\mathrm{Str\,}(\mathrm{Str\,}(f)g)=\mathrm{Str\,}(f)\mathrm{Str\,}(g)$};
											\item[(iii)] for all $\alpha\in K$, we have $\mathrm{Str\,}(\alpha g)=\alpha \mathrm{Str\,}(g)$;
											\item[(iv)] For all $f$, $g$ homogeneous elements in $K^{\mathcal{V}}\{X; Y\}$, we set \[\mathrm{Str\,}(fg)=(-1)^{\deg(f)\deg(g)}\mathrm{Str\,}(gf) {~~and ~~\mathrm{Str\,}(f)\mathrm{Str\,}(g)=(-1)^{\deg(f)\deg(g)}\mathrm{Str\,}(g)\mathrm{Str\,}(f)}.\]
										\end{enumerate}
										The elements of $STR^{\mathcal{V}}\langle X; Y\rangle$ are called pure supertrace polynomials.
									\end{definition}

									\begin{definition}We denote by ${K}^{\mathcal{V}}_{STR}\{X;Y\}$ the free supertrace graded algebra in $\mathcal{V}$ that is the supertrace {$\mathbb{Z}_2$-graded} algebra generated by $K^{\mathcal{V}}\{X;Y\}$, $STR^{\mathcal{V}}\langle X;Y\rangle$ and the symbol $r$ with the conditions:
										\begin{enumerate}
											\item[(i)] $STR^{\mathcal{V}}\langle X;Y\rangle$ is supercentral in ${K}^{\mathcal{V}}_{STR}\{X;Y\}$;
											\item[(ii)] $\mathrm{Str\,}(1)=r$ is central.
										\end{enumerate}
									\end{definition}
									
									Remark the algebras with supertrace in $\mathcal{V}$ form a variety (one includes the supertrace in the signature of the variety). 
									Unless otherwise stated we shall assume that $X$ and $Y$ are infinite countable set{s}. It is immediate that  ${K}^{\mathcal{V}}_{STR}\{X;Y\}$ contains the free algebra $K^{\mathcal{V}}\{X;Y\}$ as a subalgebra. 
									
									\begin{definition}
										The subalgebra $G_S\{X;Y\}$ of ${K}^{\mathcal{V}}_{STR}\{X;Y\}$ generated by the set 
										\[
										\{g(x_1,\dots,x_n; y_1,\dots,y_m),\mathrm{Str\,}(g(x_1,\dots,x_n; y_1,\dots,y_m))\mid g(x_1,\dots,x_n; y_1,\dots,y_m)\in K^{\mathcal{V}}\{X;Y\}\}
										\]
										is called \textsl{algebra of generalized supertrace polynomials} in the variety $\mathcal{V}$ and its elements are called \textsl{supertrace polynomials}.
									\end{definition}
									
									It is clear that $G_S\{X;Y\}$ is spanned by the generalized monomials of the form
									\begin{equation}\label{eleG}
									\widehat{a}_0\mathrm{Str\,}(a_1)\cdots \mathrm{Str\,}(a_t),
									\end{equation}
									where $t=1$, 2, \dots{}, and $a_0$, $a_1$, \dots, $a_t$ are monomials in $K^{\mathcal{V}}\{X;Y\}$.

									The supertrace monomial from (\ref{eleG}) is of \textsl{degree} $\deg a_0+\deg a_1+\dots+\deg a_t$ and the degree of a supertrace polynomial is the largest degree of a monomial that appears with non-zero coefficient in it.
									
									\begin{definition}
										Let $(A,\tau)$ be a supertrace algebra in $\mathcal{V}$. The supertrace polynomial in $G_S\{X; Y\}$ denoted by $f=f(x_1,\ldots,x_n, ; y_1,\dots,y_m)$ is a supertrace identity of $(A,\tau)$ if, replacing $\mathrm{Str\,}$ by $\tau$, we have $f(a_1,\ldots,a_n; b_1,\dots,b_m)=0$, for all {admissible substitution} $a_1$, \ldots, $a_n$ in $A^{(0)}$ and $b_1$, \ldots, $b_n$ in $A^{(1)}$. We denote by $Id_{\mathrm{Str\,}}(A,\tau)$ (or simply $Id_{\mathrm{Str\,}}(A)$) the set of all supertrace identities for $(A,\tau)$.  
									\end{definition}
									
									Clearly, all identities for $A$, considered as elements of $G_S\{X; Y\}$, are in $Id_{\mathrm{Str\,}}(A,\tau)$. Moreover, as in the ordinary case, an ideal $I$ of the algebra $G_S\{X; Y\}$ is called \textsl{$T_{\mathrm{Str\,}}$-ideal} (or \textsl{$T$-ideal with supertrace}) if, for every supertrace polynomial $f(x_1,\ldots,x_n;y_1,\ldots,y_m) \in I$ and for any {admissible substitution of } polynomial $h_1,\ldots, h_n \in {K^{\mathcal{V}}}\{X;Y\}^{(0)}$, $g_1,\ldots, g_m \in {K^{\mathcal{V}}}\{X;Y\}^{(1)}$, the supertrace polynomial $f(h_1, \dots, h_n;g_1,\ldots, g_m)$ is contained in $I$. Therefore, it is possible to verify that for each supertrace algebra $(A,\tau)$, the ideal of its supertrace identities is a $T_{\mathrm{Str\,}}$-ideal, denoted by \( \mathrm{Id}_{\mathrm{Str}}(A) \). Using a standard Vandermonde argument, it is possible to show that, in characteristic 0, any $T_{\mathrm{Str\,}}$-ideal is generated by its multilinear supertrace polynomials. 
									
									Finally, we will mention an important result from \cite{berelsuper}. Assume that $R$ is a trace algebra and consider the supertrace algebra $R':=R\otimes E$. Given a multilinear supertrace polynomial $f(x_1,\ldots,x_n,y_1,\ldots, y_m)\in G_S\{X;Y\}$ we make the formal subtitution $y_i=e_i\bar{y}_i$ under the assumptions that the $e_i$ are supercentral and the $\bar{y}_i$ are degree zero. Thus, the $e_i$ anticommute with each other and commute with everything else. Since $\mathrm{Str\,}$ is linear over the supercenter, we are allowed to write
									$$f(x_1,\ldots, x_n,e_1\bar{y}_1,\ldots, e_m\bar{y}_m)=g(x_1,\ldots, x_n,\bar{y}_1,\ldots, \bar{y}_m)e_1\cdots e_m.$$
									In this case, $g$ can be seen as a trace polynomial and we write $\theta(f)=g$.
									
									\begin{lemma}\label{berel2.2sup}\cite[Lemma 2.2]{berelsuper}
										Let $R{=R^{(0)}\oplus R^{(1)}}$ {be a supertrace algebra with}  supercenter $\mathfrak{Z}$ and assume that $ann\mathfrak{Z}^{(1)}$ is the trivial ideal. Then for all multilinear supertrace polynomials $f$, $f$ is a supertrace identity for $R$ if and only if $\theta(f)$ is a trace identity for $R^{(0)}$.
									\end{lemma}

									\begin{corollary}\cite[Corollary 2.4.]{berelsuper}
										Let $R{=R^{(0)}\oplus R^{(1)}}$ {be a supertrace $\mathbb{Z}_2$-graded algebra with} supercenter $\mathfrak{Z}$ and assume that $ann\mathfrak{Z}^{(1)}$ is trivial. Let $\{f_i\}_{i\in I}$ be a set of multilinear trace polynomials which forms a basis for the trace identities of $R^{(0)}$. Then the set $\{g\in G_S\{X;Y\}\mid \theta(g)=f_i,\mbox{ some } i\in I\}$ forms a basis for the supertrace identities of $R$.
									\end{corollary}

									\subsection{Brief review on Invariant theory}\label{invariant} Let $A^k=A\oplus \cdots\oplus A$ be the direct sum of $k$ copies of an algebra $A$ and let $\mathcal{G}\subseteq End_KA$ be an algebraic group. The group $\mathcal{G}$ acts diagonally on $A^k$. We consider the algebra $\mathcal{A}$ of polynomial functions from $A^k$ to $A$, then $\mathcal{G}$ acts on $\mathcal{A}$ as follows. If $f \in \mathcal{A}$, then 
									\begin{equation}\label{act}
									f^\mathfrak{g}(a_1, \ldots, a_k) = \mathfrak{g}\cdot f(\mathfrak{g}^{-1}\cdot a_1, \ldots, \mathfrak{g}^{-1}\cdot a_k),
									\end{equation}
									for every $\mathfrak{g}\in \mathcal{G}$ and  $a_i\in A$. {Therefore, we say $f$ is a $\mathcal{G}$-equivariant map, if  $f=f^\mathfrak{g}$ for all $\mathfrak{g}\in \mathcal{G}$. } 
									
									Let $K[A^k]$ be the algebra of polynomial functions from $A^k$ to $K$. Of course, $K[A^k]$ is a subalgebra of $\mathcal{A}$. The description of the generators of the algebra of invariants
									\[
									K[A^k]^\mathcal{G} = \{f\in K[A^k]\mid f^\mathfrak{g}=f \mbox{\textrm{ for every }}\mathfrak{g} \in \mathcal{G}\}
									\]
									is one of the main tasks of the classical Invariant Theory. It is often referred to as the First Fundamental Theorem of Invariant theory (we recall that the Second Fundamental Theorem deals with the description of the relations among the generators). The answer to such problem is known only in few cases. For instance, such description in the (ungraded) algebra of $n\times n$-matrices  is obtained in terms of trace polynomials, see \cite{procesitr}. In addition, we mention that it  is possible to describe the invariants of the $n \times n$ matrix algebra with either the symplectic or the transpose involution, and also for the Cayley algebra (see \cite{procesitr,scwars1, scwars2}). The last result is also valid over an infinite field of odd characteristic (see \cite{zubiva}). Moreover, the generic trace ring for $M_n(E)$ is a ring of invariant for $GL_n(K)$ acting on a certain free supercommutative ring, and an analogous result for the superalgebra $M_{a,b}(E)$ was also obtained by Berele in \cite{berelsuper}. The same author, in \cite{berelez2}, used the $\mathbb{Z}_2$-graded invariant theory to study the trace identities of $M_{a,b}(E)$, and, hence, a new version of the formal superstructure for $M_{a,b}(E)$ was considered. For some variations of these results, we address the interested reader to the references \cite{berelequeer,berelecolor} for more details on the topic.
									
									Let $A$ be an $n$-dimensional \textit{power associative} algebra, i.e., each subalgebra of $A$ generated by a single element is associative while $A$ is non-associative. Clearly the algebra $A\otimes F$ is power associative for any field extension $F$ of $K$. We choose $F$ to be the field of fractions of the commutative polynomial ring $K[\xi_1,\dots, \xi_n]$. Given a basis $\{a_1,\ldots, a_n\}$ of $A$, we consider the following element of $A\otimes F$, say $v:=\sum_{i=1}^n\xi_ia_i$, and we refer to $v$ as the \textsl{generic element} of $A$. Now consider the minimal polynomial of $v$
									\begin{equation}\label{pol.minimal}
									m_{v}(t)=t^s-\sigma_1(v)t^{s-1}+\cdots+(-1)^s\sigma_s(v);
									\end{equation}
									remark the coefficients $\sigma_i(v)$ are elements of  $K[\xi_1,\dots, \xi_n]$. Assuming $a=\sum_i \alpha_ia_i\in A$, we denote by $\sigma_i(a)$ the result of the substitution  $\xi_i\mapsto \alpha_i$ in $\sigma_i(v)$. The \textit{generic minimal polynomial} of $a$ will be then 
									\begin{equation}
									\label{ham_cayley}
									m_a(t)=t^s-\sigma_1(a)t^{s-1}+\dots+(-1)^s\sigma_s(a).
									\end{equation}
									It is well known that $m_a(a)=0$, see for example \cite[Chapter VI, section 3]{Jacobson2}. The \textsl{generic trace} of $a$, denoted by $\mathrm{Trd\,}(a)$, is  $\sigma_1(a)$ and the \textsl{degree} of the algebra $A$ is equal to the degree of the  polynomial (\ref{pol.minimal}). We note that $m_a(t)$ is independent of the choice of the basis of $A$ over ${K}$ and of the field extension of $K$. The generic trace is ${K}$-linear and $\mathrm{Trd\,}(1)=s$; moreover we have that the coefficients of $m_{v}(t)$ are invariant under the action of the automorphisms group of $A$.
									
									{We would like to draw the reader's attention to the fact that the generic trace may sometimes not satisfy the specific definition of a trace that we are using in this paper. However, there are many concrete algebras where the generic trace indeed satisfies our definition. For example, the generic trace on central simple associative algebras, alternative algebras, and Jordan algebras is symmetric, associative, and nondegenerate. For further details, see \cite[Lemma (32.4)]{knus1998book}.}
									
									From now on, if $A$ is a finite dimensional power associative algebra endowed with a generic trace $\mathrm{Trd\,}$. Thus, we will always consider the pair $(A,\mathrm{Trd\,})$ that is a trace algebra. We denote by $Tr_k^\prime$ the subalgebra of $K[A^k]$ generated by the elements $\mathrm{Trd\,}(p)=\mathrm{Trd\,}\circ p$ where $p\in K\{x_1,\ldots, x_k\}$, the free algebra on the set of $k$ free generators $x_1$, \ldots, $x_k$. Notice that $Tr_k^\prime\subseteq K[A^k]^\mathcal{G}$, if $\mathcal{G}\subseteq End_KA$. 
									
									The next {hypothesis} will be crucial in what follows.
									\begin{enumerate}
										\item[($\mathcal{P}_1$)] Let $\mathcal{G}$ be an algebraic group. If {${\mathcal{G}}$ is a subgroup of the group of \textcolor{blue}{graded} automorphisms of $A$, denoted by $Aut_KA$,} then $Tr_k^\prime= K[A^k]^\mathcal{G}$.
										\item[($\mathcal{P}_2$)] The generic trace on $A$ is non-degenerate, i.e., $\mathrm{Trd\,}(uv)=0$, for all $v\in A$, implies that $u=0$.
										\item[($\mathcal{P}_3$)] The group 
										${\mathcal{G}}$ is linearly reductive.
									\end{enumerate}

									{
										Assume we have a $\mathbb{Z}_2$-grading on $A$ so that $A=A^{(0)}\oplus A^{(1)}$, $\dim(A^{(0)})=a$ and $\dim(A^{(1)})=b$. Thus, we take $\mathcal{B}_A=\{u_1,\ldots,u_a,v_1,\ldots, v_b\}$ a basis of homogeneous elements of $A$ satisfying $\deg(u_i)=0$ and $\deg(v_i)=1$. Moreover, let $S=\{s_{i1},s_{i2},\ldots,s_{i(a+b)}\mid i\geq1\}$ and $T=\{t_{i1},t_{i2},\ldots,t_{i(a+b)}\mid i\geq1 \}$ be two countable sets where the elements of $S$ are central, the elements of $T$ anticommute and $[s,t]=0$, for all $s\in S$ and $t\in T$. 
										Indeed, the associative algebra generated by $S$ and $T$ is the free supercommutative algebra over $K$ that will be denoted by $\Delta_S$ from now on. Using the same notation of Corollary \ref{embsup}, we consider $\mathcal{E}$, the graded subalgebra of $A\otimes \Delta_S$ generated by of all generic elements of $A$ over $\Delta_S$ given by
										\begin{equation*}
										X_j=u_1\otimes s_{1j}+\cdots+u_{a}\otimes s_{aj}+v_1\otimes t_{(a+1)j}+\cdots+v_{b}\otimes t_{(a+b)j}
										\end{equation*}
										and
										\begin{equation*}
										Y_j=u_1\otimes t_{1j}+\cdots+u_{a}\otimes t_{aj}+v_1\otimes s_{(a+1)j}+\cdots+v_{b}\otimes s_{(a+b)j}.
										\end{equation*}}

									Following the  notation above, we will always consider the pair $(A\otimes \Delta_S,\mathrm{Strd\,})$ that is a supertrace algebra, where $\mathrm{Strd\,}=\mathrm{Trd\,}\otimes id$. The set $STR[X_1,Y_1,\ldots]$ denotes the superring generated by $\mathrm{Strd\,}(X_i)$ and $\mathrm{Strd\,}(Y_J)$, where $X_i$ and $Y_j$ are the generic elements given above. On the other hand, $TR[X_1,\bar{Y}_1,\ldots]$ is the ring generated by $\mathrm{Trd\,}(X_i)$ and $\mathrm{Trd\,}(\bar{Y}_j)$, where $X_i$ and $\bar{Y}_j$ are the generic elements in $A$. In other words, the latter are copies of $X_i$ and $Y_j$ considered without grading. Notice that $STR[X_1,Y_1,\ldots]$ is a subalgebra of the free supercommutative algebra $\Delta_S$ generated by $s_{ij}$, $t_{ij}$ whereas $TR[X_1,\bar{Y}_1,\ldots]$ is a subalgebra of the polynomial algebra in the variables $s_{ij}$, $\bar{t}_{ij}$ that will be denoted by $K[S;\bar{T}]$.
									
									{Given a homomorphic image of the free algebra $A$ with multigrading given by the homogeneous degree, we denote by $mult(A)$ the space of multilinear elements of $A$, i.e., the space spanned by homogeneous elements of degree $\leq 1$ in each variable. It turns out the map $\theta$ is a vector space isomorphism from $mult(G_S\{X;Y\})$ to $mult(G\{X;\bar{Y}\})$. Moreover, Lemma \ref{berel2.2sup} implies that $\theta$ preserves $T$-ideals. 	}
									
									{	\begin{lemma}\label{Berele223}
											The following diagram commutes:
											\begin{equation*}
											\xymatrix{mult(K\{X_1,Y_1,\ldots\}) \ar[d]_{\mathrm{Str\,}}\ar[r]^{\theta} & mult(K\{X_1,\bar{Y}_1,\ldots\})\ar[d]^{\mathrm{Tr\,}}\\
												mult(\Delta_S)\ar[r]_{\theta} & mult(K[S;\bar{T}]).\\ }
											\end{equation*}
											Moreover, the map $\theta\colon mult(STR[X_1,Y_1,\ldots])\to mult(TR[X_1,Y_1,\ldots])$ is an isomorphism of $\mathcal{G}$-module, where $\mathcal{G}=Aut_K^{\mathbb{Z}_2}(A)$ is the group of the graded automorphism of $A$.
										\end{lemma}\begin{proof}
										Firstly, the linearity of formal supertrace  $\mathrm{Str}$ on $K\{X_1,Y_1,\ldots\}$ and the definition of $\theta$ imply that 
										\begin{equation*}
										\xymatrix{mult(K\{X_1,Y_1,\ldots\}) \ar[d]_{\mathrm{Str\,}}\ar[r]^{\theta} & mult(K\{X_1,\bar{Y}_1,\ldots\})\ar[d]^{\mathrm{Tr\,}}\\
											mult(STR[X_1,Y_1,\ldots])\ar[r]_{\theta} & mult(TR[X_1,\bar{Y}_1,\ldots]).\\ }
										\end{equation*}
										is a commutative diagram. Thus, the proof of the first part follows from the extension from $\theta$ to $mult(\Delta_S) \to mult(K[S;\bar{T}])$. For this, we take $f(s^{\alpha},t^{\alpha})=f(s_1,\ldots,s_q,t_1,\ldots,t_l)\in mult(\Delta_S)$, where $s_1$, \ldots, $s_q$ are variable in $S$ whereas $t_1$, \ldots, $t_l$ lie in $T$. Let $e_\alpha=(e_1,\ldots,e_l)$ be supercentral elements in $T$ with non-zero product. Under these conditions, $\theta(f)\in mult(K[S;\bar{T}])$ can be uniquely defined so that 
										$$\theta(f(s^{\alpha},e_\alpha \bar{t}^{\alpha}))=\theta(f)(s^{\alpha}, \bar{t}^{\alpha})e_1\cdots e_l.	$$
										This conclude the first part. 
										
										Now according to Proposition \ref{actionofGsup}, the elements in $\mathcal{G}$ act on $\Delta_S$ via the representation $\xi_g=\Phi_g\circ(\zeta^{g}_A)^{-1}$. Here $\Phi_g$ act on $A$ while $\zeta^{g}$ act on  $ A \otimes \Delta_S$. Due to this comment, for each $g\in \mathcal{G}$, if $1\leq i\leq a$, then $g\cdot s_{ij}$ is the entry of the basic element $u_i$ in the generic element $g\cdot X_j$; case $a< i\leq n$, $g\cdot s_{ij}$ equals to the entry of element $v_i$ of the generic element $g\cdot Y_j$. Already if $1\leq i\leq a$, then $g\cdot t_{ij}$ is the entry of element $u_i$ of the generic element $g\cdot Y_j$. Finally if $a< i\leq n$, then $g\cdot t_{ij}$ equals to the entry of the basic element $v_i$ of the generic element $g\cdot X_j$. Following from definition of the $\mathcal{G}$-action that $g\cdot \theta(s_{ij})=\theta(g\cdot s_{ij})$ and $g\cdot \theta(t_{ij})=\theta(g\cdot t_{ij})$. Notice that every element in $\Delta_S$ is written as a linear combination of elements given by $k=s_{\alpha_1}\cdots s_{\alpha_r}t_{\beta_1}\cdots t_{\beta_l}$, with $\alpha_1\leq \cdots \leq \alpha_r$ and 
										$\beta_1<\cdots<\beta_l$. Thus, on such last element we have $\theta(k)=\theta(s_{\alpha_1})\cdots \theta(s_{\alpha_r})\theta(t_{\beta_1})\cdots \theta(t_{\beta_l})$. As the action on $g$ comes from a graded isomorphism we have
										\begin{eqnarray*}
											g\cdot\theta(k)&=&g\cdot(\theta(s_{\alpha_1})\cdots \theta(s_{\alpha_r})\theta(t_{\beta_1})\cdots \theta(t_{\beta_l}))\\
											&=&\theta(g\cdot s_{\alpha_1})\cdots \theta( g\cdot s_{\alpha_r})\theta( g\cdot t_{\beta_1})\cdots \theta(g\cdot t_{\beta_l}).
										\end{eqnarray*}
										The action of $g$ preserve the lower indices, and so this equals
										$$\theta(g\cdot s_{\alpha_1})\cdots \theta( g\cdot s_{\alpha_r})\theta( g\cdot t_{\beta_1})\cdots \theta(g\cdot t_{\beta_l})=\theta(g\cdot k).$$
										We are done.
									\end{proof}
								}

								{According to the equality $K[S,\bar{T}]^\mathcal{G}= K[A^k]^\mathcal{G}=TR[X_1,\bar{Y}_1,\ldots]$, from the previous lemma we get $mult(\Delta_S)^\mathcal{G}=\theta^{-1}(mult(K[S;\bar{T}])^\mathcal{G})$. In other words, the condition $(\mathcal{P}_1)$ is equivalent to:
									\begin{enumerate}
										\item[($\mathcal{P}_1^\prime$)] Let $\mathcal{G}$ be an algebraic group. If ${\mathcal{G}}$ is the group of the graded automorphism of $A$, then $STR[X_1,Y_1,\ldots]=\Delta_S^\mathcal{G}$.
									\end{enumerate}}

									\subsection{The Main Result}\label{sectioembeddin} Throughout this subsection $A$ will denote a power associative {graded-simple and central algebra} of finite dimension $n$ over $K$, equipped with a {$\mathbb{Z}_2$}-grading. Moreover, $\mathcal{G}$ will be the group of graded automorphisms of $A$. We assume $A$ is an algebra with generic trace, denoted by $\mathrm{Trd}$, such that $(A, \mathrm{Trd\,})$ satisfies the conditions ($\mathcal{P}_1$)-($\mathcal{P}_3$). 
									
									Our main goal is to prove that a {graded} algebra $R$, belonging to the same variety of $A\otimes E$, can be embedded into an algebra $A\otimes \Xi$, for an appropriate supercommutative associative $K$-algebra $\Xi$, if and only if $R$ satisfies all supertrace identities of the algebra $A$ (see Theorem \ref{main.theor}); such embedding is a graded homomorphism, too.
									
									Using a standard argument on extensions of the scalars we can assume, without loss of generality, that $K$ is an algebraically closed field, {that means $A$ is a graded-central-simple algebra}. For, assume the result we aim at holds for the algebraic closure $\bar{K}$ of the field $K$. Then for every $K$-algebra $R$ we consider the $\bar{K}$-algebra $R_{\bar{K}}=R\otimes \bar{K}$ equipped with the grading induced by $\mathbb{Z}_2$-grading of $R$. The latter algebra can be embedded into $A_{\bar{K}}\otimes_{\bar{K}} \Xi$ as a $\bar{K}$-algebra for some commutative-associative $\bar{K}$-algebra $F$. Hence 
									\[
									R\hookrightarrow R\otimes_K\bar{K}\hookrightarrow (A\otimes_K{\bar{K}})\otimes_{\bar{K}} \Xi\simeq A\otimes_K ({\bar{K}}\otimes_{\bar{K}}\Xi)\simeq A\otimes_K\Xi.
									\]
									Therefore $R$ can also be embedded into $A\otimes \Xi$, as a graded $K$-algebra.

									\begin{proposition}\label{prop.trace.esc.}
										Let $m$ be a graded monomial in the variables $x_1,\ldots, x_m$ of degree $0$ in the free algebra $K^{\mathcal{V}}\{X;Y\}$, and $y_1,\ldots,y_l$ of degree $1$. Suppose that $m$ is linear in one of its variables, say $y_l$. Then there exists a monomial $m^\prime$  in the free algebra $K\{x_1, \dots, x_m, y_1, \dots, y_{l-1}\}$ such that $m'$ does not depend on $y_l$ and 
										\[
										\mathrm{Str\,}(m)=\mathrm{Str\,}(m^{\prime} y_l),
										\]
										in $G_S\{X;Y\}$. The result also holds for variables of degree $0$.
									\end{proposition}
									\begin{proof}
										The proof follows word by word, with small modifications, the ideas of  \cite[Proposition 4.1]{CDP}, and we shall omit it here. 
									\end{proof}
									
									\begin{corollary}\label{tr. escolha}
										Assume that $h\in G_S\{X;Y\}$ is a pure supertrace polynomial in the variables $x_1$, \dots, $x_m$, $y_1$, \dots, $y_l$ that is linear in some of them, say $y_l$. Then there exists a polynomial $f\in G_S(x_1, \dots, x_m, y_1, \dots, y_{l-1})$ such that $
										\mathrm{Str\,}(h)= \mathrm{Str\,}(fy_l)$.
									\end{corollary}

									\begin{remark}\label{good}
										\label{rem:S=Delta}
										Assume that {$A\otimes E=(A\otimes E)^{(0)}\oplus (A\otimes E)^{(1)}$} is a {graded} algebra. Let $(A\otimes E)^{\infty}$ stand for the vector space of all sequences $\{a_i\}_{i\in I}$ of elements $a_i=a_i^{(0)}+a_i^{(1)}$(in $A\otimes E$) such that all but finitely many are 0. 
										
										We will show that the algebra of polynomial maps from $(A\otimes E)^{\infty}$ to $(A\otimes E)$ can be identified with $A\otimes {S}${, for some supercommutative-associative algebra $S$}. We denote by $X_j$ (resp. $Y_j$) the projection $\{a_i\}_{i\in I}\mapsto a^{(0)}_j$ (resp. $a^{(1)}_j$) on the homogeneous component of degree $(0)$ (resp. $(1)$) of the $j$-th entry of the sequence $\{a_i\}_{i\in I}$. We choose a basis  {$\mathcal{B}_A=\{u_1,\ldots,u_a,v_1,\ldots, v_b\}$ a basis of homogeneous elements of $A$ satisfying $\deg(u_i)=(0)$ and $\deg(v_i)=(1)$,} and let us denote by  {$\mathcal{B}_A=\{u^{\ast}_1,\ldots,u^{\ast}_a,v^{\ast}_1,\ldots, v^{\ast}_b\}$} its dual basis. We denote by $S$ the subalgebra of the space of algebra homomorphisms from $(A\otimes E)^\infty$ to $E$ generated by all projections $s_{ij} = u^{\ast}_i\circ X_j$, $s_{(a+i)j} = v^{\ast}_i\circ Y_j$, $t_{ij} = u^{\ast}_i\circ Y_j$ and $t_{(a+i)j} = v^{\ast}_i\circ X_j$. We attach to each $a\otimes s\in A\otimes S$ the map $x\mapsto s(x)a$, $x\in (A\otimes E)^\infty$ that is a \textsl{polynomial map} from $(A\otimes E)^\infty$ to $(A\otimes E)$. We identify $a\otimes s$ with its corresponding polynomial map that is $A\otimes S$ coincides with the algebra of polynomial maps from $(A\otimes E)^\infty$ to $(A\otimes E)$ with the grading induced by $E$. It is also easy to see $S\cong  S_A(K\{X;Y\})=\Delta_S$, too, using the notation of the previous section. 
										
										The group $\mathcal{G}$ of automorphisms of $A$ acts canonically on the algebra of maps $(A\otimes E)^\infty\to (A\otimes E)$ via $\mathfrak{g}\cdot f(x) = \mathfrak{g}f(\mathfrak{g}^{-1}x)$; it is clear that if $f$ is a polynomial map, then $\mathfrak{g}\cdot f$ is a polynomial map as well. Moreover $\mathcal{G}$ acts also on $A\otimes 
										\Delta_S$ via the $A$-universal supermap (see the remark before Proposition \ref{actionofGsup}). We recall that, because $A$ is unitary, we shall often identify without saying it $K$ and $K\cdot 1\subseteq A$. If $p$ is a homogeneous element of $A\otimes S$, then $\mathrm{Strd}\circ p$ is either a central or an anticommutative element in $ A\otimes S$, depending on the $\mathbb{Z}_2$-degree of $p$. This defines a supertrace on $A\otimes S$. The natural supertrace map on $A\otimes \Delta_S$ induced by the trace map of $A$ together with the graded isomorphism $S\cong \Delta_S$ yields hence a graded isomorphism of algebras with supertrace $\varphi\colon A\otimes S\to A\otimes \Delta_S$. Recall that for every $\mathfrak{g}\in \mathcal{G}$ and $f\in A\otimes S$ it holds $\varphi(\mathfrak{g}\cdot f) = \mathfrak{g}\varphi(f)$. Thus we identify $A\otimes S$ and $A\otimes \Delta_S$. 
										
										The projections 
										
										\[X_j=u_1\otimes s_{1j}+\cdots+u_{n}\otimes s_{aj}+v_{1}\otimes t_{(a+1)j}+\cdots+v_{b}\otimes t_{(a+b)j}\] and 
										\[Y_j=u_1\otimes t_{1j}+\cdots+u_{n}\otimes t_{aj}+v_{1}\otimes s_{(a+1)j}+\cdots+v_{b}\otimes s_{(a+b)j}\] are polynomial $\mathcal{G}$-invariant maps. We denote by $\mathcal{T}$ the graded subalgebra of $A\otimes \Delta_S$ generated by the equivariant maps, that is $(A\otimes \Delta_S)^\mathcal{G}$. Clearly, each $X_j, Y_j$ belongs to $\mathcal{T}$. 
									\end{remark}
									
									In the light of the previous remark, in the sequel we shall identify $\Delta_S$ with $S$. 
									
									\begin{lemma}\label{gen_mon}
										The superalgebra $\mathcal{T}$ is generated, as an algebra, by the elements $X_i$, $Y_i$ and by the {super}traces $\mathrm{Strd\,} m(X_1, \ldots,X_m,Y_1,\ldots, Y_t)$, where $m(X_1, \ldots,X_m,Y_1,\ldots, Y_t)$ runs over all monomials in the free unitary {graded} algebra $K\{X;Y\}$.
									\end{lemma}
									\begin{proof}
										Let $f\colon (A\otimes E)^k\rightarrow (A\otimes E)$ be an equivariant map. Assume, without loss of generality, that $f$ is multilinear. We consider the invariant $\bar{f}\colon (A\otimes E)^{k+1}\rightarrow E$ defined as follows:
										\[
										\bar{f}=\mathrm{Strd\,}(f\cdot Y_{t+1}).
										\]
										As $A$ satisfies the condition $(\mathcal{P}_1)$, we obtain that $\bar{f}$ lies in $STR[X_1,{Y}_1,\ldots]$, therefore it can be written as 
										\[
										\bar{f}=\mathrm{Strd\,}(f^\prime\cdot Y_{t+1}),
										\]
										for some $f^\prime\in G_S\{X;Y\}$. Now the lemma follows from condition $(\mathcal{P}_2)$.
									\end{proof}
									
									\begin{lemma}\label{free algebra}
										The {graded} algebra $\mathcal{T}$ is the free {graded} algebra with supertrace in the variety of the graded algebras with supertrace satisfying all supertrace identities of $A\otimes E$.
									\end{lemma}
									\begin{proof}
										In order to prove that $\mathcal{T}\simeq G_S\{X;Y\}/Id_{\mathrm{Str\,}}(A\otimes E)$ we first observe that $\mathcal{T}$ is a graded subalgebra of $A\otimes\Delta_S$ preserving supertraces, and Lemma \ref{berel2.2sup} implies that $Id_{\mathrm{Str\,}}(A\otimes \Delta_S)=Id_{\mathrm{Str\,}}(A\otimes E)$. Now consider the supertrace preserving graded epimorphism $\varphi$ from $G_S\{X;Y\}$ to $\mathcal{T}$ defined by $\varphi(x_i)=X_i$ and $\varphi(y_i)=Y_i$, for all $i$. Clearly
										$Id_{\mathrm{Str\,}}(A\otimes E)\subseteq \ker\varphi$. On the other hand, let $f\in \ker\varphi$, then $f(X_1,\dots,X_m,Y_1,\dots,Y_l)=0$ for the equivariant projections $X_1,\dots,X_m,Y_1,\dots,Y_l$ of $A\otimes\Delta_S$ respecting the $\mathbb{Z}_2$-grading. Therefore for each $m+l$-tuple $(a_1,\ldots,a_{m+l})$ of elements in $A\otimes E$ the equation \[0=f(X_1,\dots,X_m,Y_1,\dots,Y_l)(a_1,\ldots,a_{m+l})=f(a^{(0)}_1,\ldots,a^{(0)}_m,a^{(1)}_{m+1},\ldots, a^{(1)}_{m+l})\] holds, for any elements $a^{(0)}_1,\ldots,a^{(0)}_m$ in $(A\otimes E)^{(0)}$ and $a^{(1)}_{m+1},\ldots, a^{(1)}_{m+l}$ in $(A\otimes E)^{(1)}$. In other words, $f$ is a supertrace identity for $A\otimes E$ and $\ker\varphi=Id_{\mathrm{Str}}(A\otimes E)$.
									\end{proof}
									
									The Reynolds' operator $\pi^{\mathcal{G}}$ is the canonical ${\mathcal{G}}$-equivariant projection of an algebra $A$ onto an {isotypical} component of the trivial representation. Since ${\mathcal{G}}$ is linearly reductive, $\pi^{\mathcal{G}}$ commutes with all remaining ${\mathcal{G}}$-maps on $A$.  In particular the following equalities hold:
									\begin{eqnarray*}
										\pi^{\mathcal{G}}(ab)&=&a\pi^{\mathcal{G}}(b),\mbox{ if } a\in A^{\mathcal{G}}, b\in A\\
										\pi^{\mathcal{G}}(\mathrm{Trd\,}(a))&=& \mathrm{Trd\,}(\pi^{\mathcal{G}}(a)),\mbox{ if } a\in A.
									\end{eqnarray*}
									We refer the reader to the monograph \cite[pg. 154]{procesibookliegroup} for an
									extensive treatment of the topic. Lemma \ref{Berele223} implies that the Reynolds' operator property also holds for ${\mathcal{G}}$-maps on $A\otimes E$, i.e., the previous equalities are valid when we replace $\mathrm{Trd\,}$ by $ \mathrm{Strd}$.
									
									The following is the main result of the section.
									
									\begin{theorem}\label{main.theor}
										Let $A$ be a power associative graded algebra. Assume that $A$ is a finite-dimensional {graded-central-simple algebra} which satisfies the conditions ($\mathcal{P}_1$)-($\mathcal{P}_3$). Assume that ${B}$ is a {graded} algebra with supertrace belonging to the same variety of algebras as $A\otimes E$ and satisfying all supertrace identities of $A\otimes E$. Then there exists a supercommutative algebra  {$\Xi$} such that ${B}$ can be embedded into $A\otimes  {\Xi}$ as a graded $K$-algebra. Moreover,  $\rho({B})=(A\otimes  {\Xi})^{\mathcal{G}}$, where ${\mathcal{G}}$ {is the group of the graded automorphism of $A$}, and $( {\Xi},{\rho})$ is the $A$-universal supermap for ${B}$.
									\end{theorem}
									\begin{proof}		
										The algebra ${B}$ satisfies all supertrace identities of $A\otimes E$; hence, by Lemma \ref{free algebra} we get $\mathcal{T}/I$ is a presentation of ${B}$ for some supertrace ideal $I$ of $\mathcal{T}$. According to Remark \ref{good}, we can embed $\mathcal{T}$ in $A\otimes 
										\Delta_S$ and  $\mathcal{T}=(A\otimes \Delta_S)^\mathcal{G}$. 
										
										We assume $\mathcal{T}$ has infinitely many variables for
										a technical reason; this, of course, is not a restriction. Let $(I)$ be the ideal of $A\otimes \Delta_S$ generated by the set $I$. The induced map
										\[{B}=\mathcal{T}/I\rightarrow (A\otimes \Delta_S)/(I)\simeq A\otimes(
										\Delta_S/J)
										\]
										will be denoted by $\psi_R$, where $J$ is an ideal of $\Delta_S$ such that $A\otimes J=(I)$.
										Remark that both the induced ideal $(I)$ and $J$ are stable under the respective $\mathcal{G}$-actions, since so is $I$.
										Recall that, assuming the $\mathbb{Z}_2$-grading on $A\otimes \Delta_S$ induced by the  tensor product grading, we have that $(I)$ is a graded ideal of $A\otimes \Delta_S$. 
										It is clear that $(\Delta_S/J,\psi_{{B}})$ is an $A$-universal supermap for ${B}$. Since $\mathcal{G}$ is linearly reductive, and the projection $p_J\colon A\otimes 
										\Delta_S\to A\otimes 
										\Delta_S/J$ is $\mathcal{G}$-equivariant (see Corollary \ref{equi}), we have that $p_J$ carries the invariant subalgebra of $A\otimes \Delta_S$ onto the invariant subalgebra of $A\otimes \Delta_S/J$, that is $\mathcal{T}\to(A\otimes \Delta_S/J)^\mathcal{G}$ is an epimorphism. Furthermore $(I)$ is the kernel of $p_J$, thus we have $\psi_{{B}}({B})=(A\otimes \Delta_S/J)^\mathcal{G}$. 
										
										Therefore we only need to show $\psi_{{B}}$ is injective which in turn is equivalent to $(I)\cap\mathcal{T}=I$ by Corollary \ref{embsup}. The inclusion $I\subseteq (I)\cap\mathcal{T}$ is obvious. As it was done in \cite{procesi,CDP}, in order to prove the reverse inclusion, we use the Reynolds' operator. If $m(x_1,\dots, x_k)$ is a monomial in $G\{X\}$ (algebra of generalized polynomials in the variety of all non-associative algebra) and $x_i\mapsto a_i$ is a substitution of elements of $A$ such that $a_i\in I$ for at least one index $i$, say $a_k$, then $m(a_1,\dots, a_k)\in (I)$. We denote by $\mathfrak{M}$ the set of elements of $(I)$ obtained in this way. Then it is clear that every element of $(I)$ is a finite sum of elements of $\mathfrak{M}$. Let $a$ be a homogeneous element in $(I)\cap\mathcal{T}$, then there exists a finite collection of $m_i\in \mathfrak{M}$ having the same $\mathbb{Z}_2$-degree of $a$ such that $a=\sum m_i$. We choose a homogeneous variable $x$ of any $\mathbb{Z}_2$-degree which does not appear in any one of the $m_i$. According to Corollary \ref{tr. escolha} we obtain 
										\[
										\mathrm{str\,}(ax)= \mathrm{str\,}\bigl(\bigl(\sum_im_i\bigr)x\bigr)= \mathrm{str\,} \bigl(\sum_i(m^{\prime}_iu_i)\bigr).
										\]
										Here we use the symbol $\mathrm{str\,}$ to denote the supertrace induced by the generic trace $\mathrm{Trd\,}$. Observe that for each $i$ the element $m^{\prime}_i$ is linear in $x$ and $u_i\in I$ for each $i$. Applying the Reynolds' operator we obtain
										\[
										\mathrm{str\,}(ax)=\mathrm{str\,}\bigl(\sum_i\pi^\mathcal{G}(m^{\prime}_i)\cdot u_i\bigr).
										\]
										Now, for each $i$, the element $\pi^\mathcal{G}(m^{\prime}_i)$ is an invariant that is linear in $x$, hence, by Corollary~\ref{tr. escolha} and Lemma~\ref{gen_mon}, we have
										\[
										\pi^\mathcal{G}(m^{\prime}_i)=\sum_s \mathrm{str\,}({a}_{is}x)t_{is}+\sum_k\bar{x}_kw_{ik}.
										\]
										Here $t_{is}$, $w_{ik}\in\mathcal{T}$, $a_{is}$, $\bar{x}_k$ are monomials in $K\{X;Y\}$ and $x$ is a variable in $\bar{x}_k$. 
										Using once again Proposition \ref{tr. escolha} we have 
										\[
										\mathrm{str\,}((\bar{x}_kw_{ik})u_i)= \mathrm{str\,}(m_{ik}x),
										\] 
										for some $m_{ik}\in I$. Therefore we have
										\begin{eqnarray*}
											\mathrm{str\,}(ax)&=& \mathrm{str\,}\bigl(\sum_i\bigl(\sum_s (\mathrm{str\,}({a}_{is}x)t_{is}+ \sum_k\bar{x}_kw_{ik}\bigr)u_i\bigr)\\
											&=&\mathrm{str\,}\left(\left(\sum_i\left(\sum_s \mathrm{str\,}( t_{is}u_i){a}_{is}+\sum_km_{ik}\right)\right) x\right).
										\end{eqnarray*}
										By condition ($\mathcal{P}_2$), the symmetric bilinear form $\mathrm{tr}(xy)$ is non-degenerate; because $ann E^{(1)}=0$, we obtain $a=\sum_{is}(\mathrm{str\,}(t_{is}u_i)a_{is})+\sum_{ik}m_{ik}$. Since $I$ is closed under the {super}trace operation, we get $a\in I$ and we are done.
									\end{proof}

									{Now let $B$ be an algebra with supertrace in a same variety $\mathcal{V}$ of $A\otimes E$. Thus, repeating verbatim the discussion done in Section \ref{universal} and considering the free supertrace $\mathbb{Z}_2$-graded algebra $K^\mathcal{V}_{STR}\{X;Y\}$, we have that there exists an $A$-universal supermap $(S_A(B),\varrho_B)$ such that the homomorphism $\varrho\colon B\to A\otimes S_A(B)$ is compatible with supertraces, where $A$ is equipped with its generic trace. Consequently a representable functor, an action of $\mathcal{G}=Aut_K^{\mathbb{Z}_2}(A)$ on $S_A(B)$,  and an universal (supertrace preseving) supmap $\varrho_B$ whose its image is $(A\otimes S_A(B))^{\mathcal{G}}$ are obtained. We will use the notation the $A$-universal supermap for $B$ by $(S_A(B),\varrho_B)$ to emphasize that now we work with algebras with supertrace.}
									
									{We see the map $\varrho_B\colon B\to (A\otimes S_A(B))^{\mathcal{G}}$ is always surjective, and its kernel contains the elements formed by the evaluations of {$\mathrm{Id}_{\mathrm{Str}}(A\otimes E)$} that is denoted by $I_A(B)$. The latter claim implies that $I_A(B)\subseteq \ker\varrho_B$. On the other hand, Corollary \ref{mergsup} implies that $\ker\varrho_B \subseteq I_A(B)$. Summarizing we have the following consequence of the result above. 
										\begin{corollary}
											The universal map $\varrho_B\colon B\to  (A\otimes S_A(B))^{\mathcal{G}}$ of an algebra with trace is surjective, and its kernel is the ideal $I_A(B)$ (closed under trace) formed by all the evaluations of elements in $\mathrm{Id}_{\mathrm{Str}}(A\otimes E)$.
										\end{corollary}}	
										
										Observe that in the previous results we are not requiring $B$ has an unit, or if there is, we are not requiring that $\varrho_B(1_B)$ to be an unit of $A\otimes S_A(B)$. Assume that the algebra $A$ has degree $n$. Thus, if $\varrho_B(1_B)$ to be an unit of $A\otimes S_A(B)$ is non-zero, then the we have a homomorphism from the subring $\mathbb{Q}[\mathrm{Strd\,}(1_B)]=\mathbb{Q}[\mathrm{Trd\,}(1_A)]$ to $\mathbb{Q}$ mapping $\mathrm{Trd\,}(1_A)$ to $n$. Here $\mathrm{Trd\,}$ is the generic trace on $A$ and $\mathrm{Srd\,}=\mathrm{Trd\,}\otimes id$, where $id$ is the identity map on {$S_A(B)$}.
										
										\begin{remark}
											In the hypotheses of Theorem \ref{main.theor}, if $B$ is a finitely generated algebra generated by $k$ elements, then {the free supercommutative algebra $\Delta_S=K[S]\otimes E$} can be chosen so that
											$K[S]=K[s_1,\ldots,s_{k\dim(A)}]$  is finitely generated by $k\dim(A)$ elements too. Hence we get \[\text{\rm GK}((\Delta_S/J))^{(0)}\leq \text{\rm GK}((\Delta_S))^{(0)}=k\dim(A).\]
											
										\end{remark}
										
										In the light of the previous remark, because of Corollary \ref{GK2}, we immediately get the next.
										
										\begin{proposition}
											In the hypotheses of Theorem \ref{main.theor}, if $B$ is a finitely generated algebra generated by $k$ elements, then \[\text{\rm GK}({B})\leq k\dim(A).\]
										\end{proposition}
										
										If we specialize $B$ with the relatively free associative algebra of $A\otimes E$ generated by $k$ elements, it is well known (see for instance \cite{berelegk} or \cite{centroneGK}) $\text{\rm GK}({B})=ek-a$, where $e,a$ are positive integers. Hence, for sufficiently big $k$, we get $e\leq \dim(A)$. A similar result was obtained by one of the authors in \cite{centroneGK} for algebras graded by a finite group. More precisely, in the case of graded algebras, $e=\text{\rm exp}^G(B)$, where $B$ is a finite-dimensional algebra such that $Id^G(A\otimes E)=Id^G(E(B))$, and $E(B)$ is the so-called Grassmann envelope of $B$ in the light of the Representability Theorem for $G$-graded algebras (see \cite{AljadeffBelov2010}, \cite{Sviridova}, \cite{Kemer1987} and \cite{kemer}); $\text{\rm exp}^G(B)$ denotes the $G$-graded exponent of $B$.
										
										We can also relate the Gelfand-Kirillov dimension of $B$ with the even (Krull) dimension of $\SSpec(S_A(B))$, as we can see by the following result, which is an easy consequence of Corollary \ref{GK2} and the fact that the Gelfand-Kirillov dimension of a finite generated commutative associative algebra coincides with its Krull dimension. 
										
										\begin{corollary}\label{lastcor}
											Assume that $A$ is a finite-dimensional power associative {graded-central-simple algebra} which satisfies conditions ($\mathcal{P}_1$)-($\mathcal{P}_4$). Fix also $B$ a graded algebra with a supertrace belonging to the same variety of algebras as $A\otimes E$ and satisfying all supertrace identities of $A\otimes E$. Then 
											
											$$\text{\rm GK$(B)$} \leq \dim \SSpec(S_A(B)).$$
										\end{corollary}
										
										\section{Formal Smoothness for central simple algebras {over supercommutative ring}}\label{section.suave}
										
										Herein we aim to study a geometrical property related to the algebras already seen since now: the \textit{formal smoothness}. More precisely we will relate the formal smoothness of a given algebra with that of its embedding algebra.
										
										\
										
										Following \cite{procesi}, we have the next definition.
										\begin{definition}\label{S)}
											\cite[Definition 5.1]{CCGIsrael} Let $\mathcal{V}$ be a variety of $G$-graded algebras. We say ${B}$ in $\mathcal{V}$ is formally smooth if for any algebra $U$ in $\mathcal{V}$, any graded nilpotent ideal $N$ in $U$, we can lift any $G$-graded homomorphism  $\varphi\colon {B} \to U/N$ to a $G$-graded homomorphism $\varphi^\prime \colon {B} \to U$.
										\end{definition}
										
										{Handling lifting arguments for $G$-graded algebras gives rise to the problem of dealing with all ideals instead of just the graded ones. The main result of this section addresses this problem properly. }
										
										In the sequel we shall use the notation adopted before. Let $A$ be a finite-dimensional {graded-central-simple algebra} in the variety of graded algebras $\mathcal{V}$ satisfying $(\mathcal{P}_1)-(\mathcal{P}_3)$; let $\mathrm{Strd}$ be the supertrace on $A\otimes \Delta_S$ induced by the generic trace on $A$. {As it was considering in the previous sections, here $\Delta_S$ denotes the free supercomutative algebra freely generated by two countable sets, $S$ and $T$, over $K$.} Due to Lemma \ref{free algebra}, for every  graded algebra with supertrace ${B}$ in $\mathcal{V}$ satisfying all supertrace identities of $(A\otimes {E},\mathrm{Strd\,})$, there is a graded ideal $I_{{B}}$ such that ${B} = \mathcal{T}/I_{{B}}.$
										Thus, there exists an ideal $I_A(B)$ of $\Delta_S$ such that $(I_{{B}}) = A \otimes I_A(B)$, and ${S_A(B)}=\Delta_S/I_A(B)$ is the supercomutative ring in the $A$-universal supermaps for $B$. As already mentioned in the proof of Theorem \ref{main.theor}, $I_A(B)$ is a $\mathcal{G}$-stable graded ideal of $\Delta_S$.
										
										Now, assume $A$ satisfies the following additional condition:
										\begin{enumerate}
											\item[($\mathcal{P}_4$)] Any monomial inside a supertrace in $\mathcal{T}$ has degree bounded by a certain $c>0$.
										\end{enumerate}
										
										We have the following technical result.
										
										\begin{lemma}
											For all $k$, we have $(I_{{B}})^{kc} \cap \mathcal{T} \subseteq I_{{B}}^k.$
										\end{lemma} 
										\begin{proof}
											Let $\mathcal{T}_{n}^{\bf g}$ be the supertrace algebra on the generic superelements $\{X_1, \ldots, X_a,Y_1, \ldots, Y_b\}$, while $\mathcal{T}_{n+l}^{\bf h}$ denotes the supertrace algebra generated by $\{X_1, \ldots, X_a,Y_1, \ldots, Y_b,Z_1, \ldots, Z_l\}.$ Here $a+b=n$, ${\bf g}=(\underbrace{0,\ldots, 0}_{ \times a},\underbrace{1,\ldots, 1}_{ \times b})$, and $Z_1, \ldots, Z_l$ denote new $\mathcal{G}$-invariant maps of type $X_i$ or $Y_j$. 
											
											Fix a set $\{U_1,\ldots,U_l\}$ of elements in $\mathcal{T}_n^{\bf g}$ so that for each $j$ the element $Z_j$ has the same $\mathbb{Z}_2$-degree of $U_j$; we define the map $u\colon \mathcal{T}_{n+l}^{\bf h} \to \mathcal{T}_n^{\bf g}$ sending $Z_j$ to $U_j$. Notice that $u$ preserves the supertrace and induces a $\mathcal{G}$-equivariant map of the corresponding supercommutative--associative algebras $\Xi_1$ and $\Xi_2$ coming from the $A$-universal supermap for the free superalgebras $K\{x_1, \ldots, x_a,y_1, \ldots, y_b\}$ and $K\{x_1, \ldots, x_a,y_1, \ldots, y_b,z_1, \ldots, z_l\}$, respectively, where $\deg(z_j)=\deg(Z_j)$ for each $j$. Thus, the next diagram is commutative.
											
											\[\xymatrix{ A \otimes \Xi_2 \ar[r] \ar[d]_{\pi^\mathcal{G}} & A \otimes \Xi_1 \ar[d]^{\pi^\mathcal{G}} \\
												\mathcal{T}_{n+l}^{\bf h} \ar[r]^{u} & \mathcal{T}_{n}^{\bf g}} \]
											We get that any monomial 
											$$m(x_1, \ldots, x_a,y_1, \ldots, y_b,z_1, \ldots, z_l)$$ 
											in the algebra of generalized supertrace polynomials $G_S\{x_1, \ldots, x_a,y_1, \ldots, y_b,z_1, \ldots, z_l\}$ in $\mathcal{V}$. Thus if $h_1,\ldots, h_a\in  (A \otimes \Xi_1)^{(0)}$ and $k_{1},\ldots, k_b\in (A \otimes \Xi_1)^{(1)}$, we can compute 
											$$\pi^\mathcal{G}(m(h_1,\ldots, h_a, k_{1},\ldots, k_b, U_1, U_2, \ldots, U_l))$$ {in the following way: }first, we compute  $$\mathfrak{R}=\pi^\mathcal{G}(m(h_1,\ldots, h_a, k_{1},\ldots, k_b, Z_1, Z_2, \ldots, Z_l)),$$ then we substitute $Z_i$ by $U_i$.
											
											By $(\mathcal{P}_1)$ and $(\mathcal{P}_2)$, we can write $\mathfrak{R}$ in terms of the sum of invariants as
											\begin{equation}\label{mmbo}
											\mathrm{str\,}(m_1)\cdots \mathrm{str\,}(m_{t}) m_{t+1},
											\end{equation}
											and  $(\mathcal{P}_4)$ implies that the monomials $m_1$, \ldots, $m_{t}$ have degree $\leq c$. Thus, for each $1\leq i\leq {t}$, $m_i$ contains at most $c$ of the $Z$'s, hence the monomial $m_{t+1}$ contains at least $l-vc$ of the $Z$'s, where $v \leq t$ is the number of $m_i$, with $i \leq t$, containing at least one $Z$.
											
											Now, assume that the evaluation of $Z$'s by elements $U$'s lies in the ideal $I_{{B}}$. Since $I_{{B}}$ is a supertrace ideal of $\mathcal{T}$, i.e., an ideal closed under taking supertraces and that preserves its $\mathbb{Z}_2$-grading on the supertrace map, then it follows that \eqref{mmbo} lies in $I_{{B}}^{v+(l-vc)}\subseteq I_{{B}}^k$, if we take $l = kc$ as $v+|k-v|c \geq k$. But this finishes the proof of the required inclusion when the $U$'s are homogeneous elements in the grading.

											The case in which $U$'s are not all homogeneous, we consider each homogeneous summand and repeat the process described above. Thus, each summand will be in the required space $I_{R}^k$.  Therefore, each summand will be in the required space $I_{{B}}^k$, and we are done, since $I_{{B}}$ is an ideal of $\mathcal{T}$.
										\end{proof}
										
										Let $({B},\tau)$ be a supertrace algebra in $\mathcal{V}$ such that $\mathrm{Id}_{\mathrm{Str}}(A\otimes E,\mathrm{Strd\,})$ is contained in $Id_{\mathrm{Str\,}}({B},\tau)$, and a graded ideal $I$ of ${B}$ closed under taking supertraces. From the previous lemma, we get, for every power $k$, 
										\begin{equation}\label{rema1}
										(I)^{kc} \cap {B} \subseteq I^k,
										\end{equation}
										where $(S_A(B), \varrho_B)$ is the $A$-universal supermap for ${B}$, and $(I)$ denotes the ideal in $A\otimes {S_A(B)}$ generated by set $I$.
										The following theorem is the main result of this section.
										\begin{theorem}\label{mainSmooth}
											Let $A$ be a power associative graded algebra such that $A$ is a finite-dimensional {graded-central-simple algebra}  which satisfies the conditions ($\mathcal{P}_1$)-($\mathcal{P}_4$). Assume that ${B}$ is a graded algebra with supertrace belonging to the same variety of algebras as $A\otimes E$ and satisfying all supertrace identities of $A\otimes E$. If $(S_A(B), \varrho_B)$ is the $A$-universal supermap for $B$, then $B$ is formally smooth (in $\mathcal{V}$) if and only if $S_A(B)$ is formally smooth as a superring.
										\end{theorem}
										\begin{proof}
											Suppose that $B$ is formally smooth. For any  supercommutative-associative algebra $C$ and graded nilpotent ideal $N$ of $C$, we get $\gamma\colon {S_A(B)} \to C/N$ a homomorphism of supercommutative-associative algebras, and we consider the composition
											\[\xymatrix{{B}
												\ar[r]^{\varrho_B} & A \otimes {S_A(B)} \ar[r]^{\gamma_A} & A \otimes C/N} \]
											where $\gamma_A$ is the map induced by $\gamma$,  we must lift $\gamma$. We know that $A \otimes (C/N) \simeq \dfrac{A \otimes C}{A \otimes N}$, where $A \otimes N$ is a graded nilpotent ideal of $A \otimes C$. As $B$ is formally smooth, the map $f := \gamma_A\circ \varrho_B$  can be lifted to $f'\colon B \to A \otimes C$. Thus, we obtain the following commutative diagram
											\[\xymatrix{ {B} \ar[r]^{f'} \ar[d]_{\varrho_B} & A \otimes C \ar[d]^{\pi_N} \\
												A \otimes {S_A(B)} \ar[r]_{\gamma_A} & A \otimes (C/N)}. \]  
											Due to Theorem \ref{Univsup}, there is a graded homomorphism $\eta_A\colon A \otimes {S_A(B)}\to A \otimes C$ induced by $\eta\colon {S_A(B)} \rightarrow C$, such that $\eta_A \varrho_B = f'.$ Thus, $\eta$ is the lift of $\gamma$ that we wanted.

											Conversely, assume ${S_A(B)}$ is formally smooth, then we take $U \in \mathcal{V}$ a supertrace algebra such that $Id_{\mathrm{Str\,}}(A{\otimes E})\subseteq Id_{\mathrm{Str\,}}(U)$ and a graded nilpotent ideal $N$ of $U$. We will prove that each {graded} homomorphism $f\colon {B} \to U/N$ can be lifted to a {graded} homomorphism $f'\colon {B} \rightarrow U$. Let $(S_A(U), \varrho_U)$ be the $A$-universal supermap for $U$, we obtain that the diagram	
											\[\xymatrix{ & U \ar[r]^{{\varrho_{U}}} \ar[d]_{p_N} & A \otimes {S_A(U)} \ar[d]^{p_J}\\
												{B} \ar[r]^f \ar[drr]_{{\varrho_{B}}} & U/N \ar[r] & A \otimes ({S_A(U)}/J) \\  &  & A \otimes {S_A(B)} \ar@{-->}^{}[u]<1mm>} \]
											is commutative. Here,  $J$ is the graded $\mathcal{G}$-invariant ideal of ${S_A(U)}$ such that $A \otimes J = (N)$. By the universal property of ${\varrho_{B}}$, we obtain a map $\alpha\colon {S_A(B)} \rightarrow {S_A(U)}/J$. However, $J$ is usually not nilpotent, whence we cannot apply directly the hypothesis on ${S_A(B)}$. As $N$ is a nilpotent ideal of $U$, there is some $s$ such that $N^s = 0$. By \eqref{rema1},
											$$(A \otimes J)^{sc} \cap U \subseteq N^s = 0$$
											so that $U$ embeds in $A \otimes ({S_A(U)}/J^{sc})$ that preserves grading and supertrace. Now, we observe the following diagram is commutative:
											\[\xymatrix{ & U \ar[r]^{\iota} \ar[d]^{p_N} & A \otimes ({S_A(U)}/J^{sc}) \ar[d]\\
												B \ar@{-->}^{}[ur]<1mm> \ar[r]^f \ar[drr]_{{\varrho_{B}}} & U/N \ar[r] & A \otimes ({S_A(U)}/J) \\ & &A \otimes {S_A(B)} 
												\ar[u]_{\alpha_A}\ar@/_1.5cm/[uu]_{{\alpha}^{\prime}_A}  }. \] 
											Notice that we can lift $\alpha$ to  $\alpha^\prime\colon  {S_A(B)} \to {S_A(U)}/J^{sc}$. This in turn produces a supertrace preserving the {graded}  homomorphism $$\beta\colon {B} \to A \otimes ({S_A(U)}/J^{sc})$$
											given by $\beta={\alpha}^{\prime}_A\circ {\varrho_{B}}$. In addition, the image of $\beta$ is contained in the algebra of $\mathcal{G}$-invariants of $A \otimes ({S_A(U)}/J^{sc})$. As $U \hookrightarrow A \otimes ({S_A(U)}/J^{sc})$, and by surjectivity of invariants under surjective maps, the algebra of $\mathcal{G}$-invariants of $A \otimes ({S_A(U)}/J^{sc})$ is equal to $U$, giving the required lift $\beta$.
										\end{proof}
										Our next goal is to relate the formal smoothness of an algebra $B$ as in the previous theorem, with its representing superscheme. To this end, we recall some definitions as they were presented in \cite{BHP}.
										
										\begin{definition}\label{SS}\cite[Definition A.18.]{BHP}
											A morphism of affine superschemes $f\colon X\to Y$ is comm-formally smooth if for
											every affine $Y$-superscheme $\SSpec(A)$ and every nilpotent ideal $N \subseteq A$, any $Y$-morphism
											$\SSpec (A/N) \to X$ extends to an $Y$-morphism 	$\SSpec (A) \to X$. When $X$ is a superscheme over a field $\mathbb{K}$, we will say that $X$ is smooth if, and only if $f : X \to Spec({K})$ is smooth. 
										\end{definition}
										
										We now obtain the equivalence we were seeking for.
										
										\begin{corollary}\label{explicitcriteria}
											Assume that $A$ is a finite-dimensional power associative {graded-central-simple algebra}  which satisfies the conditions ($\mathcal{P}_1$)-($\mathcal{P}_4$). Let $B$ be a graded algebra with supertrace belonging to the same variety of algebras as $A\otimes E$ and satisfying all supertrace identities of $A\otimes E$. Then $\SSpec(S_A(B))$ is smooth as a superscheme if and only if $B$ is formally smooth.
										\end{corollary}
										\begin{proof}
											It is clear from the Definitions \ref{S)} and \ref{SS} that an affine superscheme $X = \SSpec(S)$ over a field ${K}$ is smooth if, and only if, $S$ is formally smooth as super commutative ${K}$-algebra.  Now, by Theorem \ref{mainSmooth} we see that $S_A(B)$ is formally smooth if, and only if, $B$ is formally smooth, from what follows our result. 
										\end{proof}
										
										{The above result helps us to provide concrete examples of algebras that are, and that are not, formally smooth, by relating them with their geometrical counterpart. For instance, Consider $A = \mathbb{O}$ the octonion algebra, the supercommutative rings $S_1 = \mathbb{C}[x_1,x_2,x_3, \theta_1 ,\theta_2]/(x_1^2-x_2x_3)$ and $S_2 = \mathbb{C}[x_1,x_2, \theta_1 ,\theta_2]/(x_2^3-x_1^2)$, {where $x_i$ and $\theta_j$ are even and odd variables, respectively,} and finally the group $G_2$ that is the automorphism group of $\mathbb{O}$. By Theorem  \ref{main.theor} and Corollary \ref{explicitcriteria}, the algebra of invariants $B_1 = (\mathbb{O} \otimes S_1)^{G_2}$ is formally smooth, while  $B_2= (\mathbb{O} \otimes S_2)^{G_2}$ is not because, by the Jacobian criteria for superschemes (\cite[Proposition A.20]{BHP}), $\SSpec(S_1)$ is smooth (actually it is a non-degenerated conic) while $\SSpec(S_2)$ is not (it is a cuspidal cubic).}
										
										\section{Acknowledgements}
										The second author wants to thank GNSAGA (Gruppo Nazionale per le Strutture Algebriche, Geometriche e le loro Applicazioni).
										
										\section{Statements and Declarations}
										
										The authors declare they have no conflicts of interests.

										\color{black}{}

									\end{document}